\documentclass[12pt]{article}
\usepackage{amsthm, amssymb, stmaryrd, amsmath, fullpage}
\usepackage[all]{xy}

\numberwithin{equation}{subsection}
\newtheorem{theorem}{Theorem}[subsection]
\newtheorem{lemma}[theorem]{Lemma}

\newtheorem{cor}[theorem]{Corollary}
\newtheorem{prop}[theorem]{Proposition}

\theoremstyle{definition}
\newtheorem{defn}[theorem]{Definition}
\newtheorem{example}[theorem]{Example}
\newtheorem{remark}[theorem]{Remark}

\newtheorem{convention}[theorem]{Convention}
\newtheorem{hypo}[theorem]{Hypothesis}

\def\AAA{\mathbb{A}}
\def\CC{\mathbb{C}}
\def\FF{\mathbb{F}}

\def\QQ{\mathbb{Q}}
\def\RR{\mathbb{R}}
\def\ZZ{\mathbb{Z}}

\newcommand{\calE}{\mathcal{E}}
\newcommand{\calF}{\mathcal{F}}
\newcommand{\calG}{\mathcal{G}}

\newcommand{\calM}{\mathcal{M}}
\def\calO{\mathcal{O}}
\def\calP{\mathcal{P}}

\def\calR{\mathcal{R}}

\newcommand{\calU}{\mathcal{U}}

\newcommand{\gothm}{\mathfrak{m}}
\newcommand{\gotho}{\mathfrak{o}}

\def\alg{\mathrm{alg}}

\def\sep{\mathrm{sep}}

\def\beq{\begin{equation}}
\def\eeq{\end{equation}}

\def\dual{\vee}

\def\del{\partial}

\def\be{\mathbf{e}}

\def\bv{\mathbf{v}}
\def\bw{\mathbf{w}}

\DeclareMathOperator{\Aut}{Aut}

\DeclareMathOperator{\charac}{char}
\DeclareMathOperator{\coker}{coker}
\DeclareMathOperator{\conv}{conv}
\DeclareMathOperator{\Ext}{Ext}

\DeclareMathOperator{\Gal}{Gal}

\DeclareMathOperator{\GL}{GL}
\DeclareMathOperator{\Hom}{Hom}
\DeclareMathOperator{\id}{id}
\DeclareMathOperator{\image}{im}
\DeclareMathOperator{\Ind}{Ind}
\DeclareMathOperator{\inte}{int}

\DeclareMathOperator{\op}{op}
\DeclareMathOperator{\qc}{qc}
\DeclareMathOperator{\qu}{qu}
\DeclareMathOperator{\rank}{rank}
\DeclareMathOperator{\Res}{Res}

\DeclareMathOperator{\Spec}{Spec}
\DeclareMathOperator{\spec}{sp}

\DeclareMathOperator{\unip}{unip}
\DeclareMathOperator{\unr}{unr}

\newcounter{fixmectr}

\begin{document}

\title{Local monodromy of $p$-adic differential equations: an overview}
\author{Kiran S. Kedlaya \\ Department of Mathematics \\ Massachusetts
Institute of Technology \\ 77 Massachusetts Avenue \\
Cambridge, MA 02139 \\
\texttt{kedlaya@math.mit.edu}}
\date{version of May 22, 2005}

\maketitle

\begin{abstract}
This primarily expository article collects together some facts
from the literature about the monodromy of differential equations
on a $p$-adic (rigid analytic) annulus, though often with simpler
proofs. These include Matsuda's
classification of quasi-unipotent $\nabla$-modules,
the Christol-Mebkhout construction of the ramification filtration,
and the Christol-Dwork Frobenius antecedent theorem. We also
briefly discuss the $p$-adic local monodromy theorem without proof.
\end{abstract}

\tableofcontents

\section{Introduction}

This 
paper is the result of an attempt to collect in one place, and present in
a uniform fashion, some disparate
results about the local monodromy of $p$-adic differential equations.
It was initiated as part of a project to establish ``semistable reduction''
for overconvergent $F$-isocrystals \cite{me-part1}; 
however, we have decided to separate
the paper from this project, as it may have independent interest.
This interest would arise from the fact that
the local monodromy of $p$-adic differential equations
intervenes both in the study of $p$-adic (rigid) cohomology, largely via
the work of Crew (e.g., see \cite{crewfin}), and in $p$-adic Hodge theory,
largely 
via the work of Berger (e.g., see \cite{berger-cst}, \cite{berger-weak}).

The purpose of this paper is mainly expository: it is intended to provide
an easy entry point into the literature on $p$-adic differential equations
for the reader familiar with rigid analytic geometry on a rather basic
level (for the most part, the only spaces being considered are annuli).
Results exposed here include:
\begin{itemize}
\item a classification of quasi-unipotent
modules with connection (due to Matsuda); 
\item
a relationship between ramification in the monodromy representation
and generic radii of convergence (based on work of Christol-Mebkhout,
Crew, Matsuda, Tsuzuki);
\item
the Frobenius antecedent theorem (due to Christol-Dwork);
\item
without proof, the $p$-adic local monodromy theorem (due to
Andr\'e, Mebkhout, and the present author).
\end{itemize}

In the remainder of this introduction, we explain a bit about what
it means for a $p$-adic differential equation to have meaningful
monodromy, then outline the structure of the paper.

\subsection{$p$-adic differential equations and their monodromy}

Let $K$ be a field of characteristic zero
complete with respect to a nonarchimedean absolute
value, whose residue field $k$ has characteristic $p>0$. Throughout
this paper, we will be considering the following situation.
We are given a rigid analytic annulus over $K$ and a ``differential
equation'' on the annulus, i.e., a module equipped with a connection
(which is automatically integrable because we are in a one-dimensional
setting). We now wish to define the ``monodromy around the puncture''
of this connection, despite not having recourse to the analytic
continuation we would use in the analogous classical setting.

We can define a monodromy representation associated to a connection if 
we can find enough horizontal sections ``somewhere''. We will only
be looking for horizontal sections on certain \'etale covers of the
annulus (what we call ``formally \'etale covers''). Specifically,
we will consider connections which on some such cover become unipotent
(filtered by submodules with trivial successive quotients), and
define a monodromy representation for these; this will give
an equivalence of categories between such ``quasi-unipotent''
modules with connection and a certain representation category.
(Beware that if the field $k$ is not algebraically closed, these
representations will only be semilinear.)

In order for such an equivalence to be useful, we need to be able
to establish conditions under which a module with connection is
forced to be quasi-unipotent. This has been done in case $K$ is
discretely valued, by the $p$-adic local monodromy theorem ($p$LMT) 
of Andr\'e \cite{andre}, Mebkhout \cite{mebkhout}
and this author \cite{me-local}. The sufficient condition in this
theorem is a so-called ``Frobenius structure'' on the connection;
like its complex-analytic analogue (variation of Hodge structure),
this extra structure arises naturally in geometric settings, and can
also be found (following examples of Dwork) in settings where the geometric
origin is a little less clear.

\subsection{Structure of the paper}

We conclude this introduction with a rundown of the contents
of the various chapters of the paper.

In Chapter~\ref{sec:ramif1}, we recall a bit of the theory
of (one-dimensional) local fields,
mostly but not exclusively in the classical case of perfect residue field.
This will be needed later to talk about monodromy representations.

In Chapter~\ref{sec:rigid}, we introduce rigid annuli,
and verify (after Gruson) that over a spherically
complete coefficient field, any coherent locally free sheaf on a
one-dimensional rigid annulus is freely generated by local sections.
This will provide a crucial link back to the literature on
$p$-adic differential equations, which is mostly phrased in terms
of modules over rings of suitably convergent power series.
We also describe the ``formally \'etale covers'' of rigid annuli
that will intervene in the study of monodromy of $\nabla$-modules.

In Chapter~\ref{sec:mono}, we define quasi-constant and quasi-unipotent
modules with connection on a one-dimensional rigid annulus,
and construct monodromy representations corresponding to such objects.

In Chapter~\ref{sec:cm}, we explain how the ramification of the monodromy
of a quasi-unipotent connection is controlled by certain ``radius of
convergence'' data; this reprises results of Christol-Dwork,
Christol-Mebkhout, Matsuda, Crew, and Tsuzuki.

In Chapter~\ref{sec:frob}, we introduce Frobenius structures,
state the $p$-adic local monodromy theorem (roughly, every
module with connection admitting a compatible Frobenius structure
is quasi-unipotent), and verify the Frobenius antecedent theorem of
Christol-Dwork.

\section{Ramification in one dimension}
\label{sec:ramif1}

We start with a quick review of ramification theory for local fields,
mostly following Serre \cite{serre}.

\setcounter{equation}{0}
\begin{convention}
When speaking of a ``discretely valued field'', we insist that the valuation
be nontrivial.
\end{convention}

\subsection{Ramification filtrations}

We recall some definitions and results from \cite[Chapter~IV]{serre}.

\begin{defn}
For $F$ a complete discretely valued field,
let $\gotho_F$ be the ring of integers of $F$, let
$\gothm_F$ be the maximal ideal of $\gotho_F$, let
$\overline{F} = \gotho_F/\gothm_F$ be the residue field of $F$,
and let $v_F: F^* \to \ZZ$ be the valuation on $F$.
\end{defn}

\begin{hypo} \label{H:comp}
For the remainder of this section, 
let $F$ be a complete discretely valued field such that
$\charac(\overline{F}) = p > 0$,
and let $E/F$ be a separable algebraic extension. Then $E$ admits a unique
valuation extending the valuation on $F$; moreover, if
$E/F$ is finite, then $E$ is complete for its valuation
\cite[Proposition~II.3]{serre}.
\end{hypo}

\begin{defn} \label{D:lower}
Assume that $E/F$ is finite Galois and that $\overline{E}/\overline{F}$
is separable (hence also Galois).
For $i \geq -1$, let $G_i$ be the subgroup of $G = \Gal(E/F)$ consisting
of those $g$ for which $v_E(a^g - a) \geq i+1$ for all $a \in \gotho_E$;
the decreasing filtration $\{G_i\}$ is called 
the \emph{lower numbering filtration} of $G$ \cite[\S IV.1]{serre}
It can be shown \cite[\S IV.2]{serre}
that $G_{-1}/G_0 \cong \Gal(\overline{E}/\overline{F})$,
that $G_0/G_1$ is cyclic of order prime to $p$, and that
$G_i/G_{i+1}$ is an elementary abelian $p$-group for $i \geq 1$.
\end{defn}

\begin{defn} \label{D:upper}
With notation as in Definition~\ref{D:lower}, 
define the function
\[
\phi_{E/F}(u) = \int_{0}^u \frac{dt}{[G_0:G_t]}.
\]
Then $\phi_{E/F}$ is a homeomorphism of $[-1, \infty)$ with itself;
let $\psi_{E/F}$ denote the inverse function.
Define the \emph{upper
numbering filtration} of $G$ by $G^i = G_{\psi_{E/F}(i)}$
\cite[\S IV.3]{serre}; its key property (``Herbrand's theorem'')
is that it commutes with formation of quotients, in that if 
$E'/F$ is a Galois subextension of $E/F$ with $H = \Gal(E'/F)$,
then the image of $G^i$ in $H$ is precisely $H^i$
\cite[Proposition~IV.14]{serre}.
\end{defn}

\begin{remark} \label{R:trans}
With notation as in Definition~\ref{D:upper},
the functions $\phi_{E/F}$ and $\psi_{E/F}$ have the following
transitivity property \cite[Proposition~IV.15]{serre}: for
$E'/F$ a Galois subextension of $E/F$,
\[
\phi_{E/F} = \phi_{E'/F} \circ \phi_{E/E'} \qquad \mbox{and}
\qquad \psi_{E/F} = \psi_{E/E'} \circ \psi_{E'/F}.
\]
\end{remark}

\begin{defn} \label{D:upper2}
Let $E$ be a Galois extension of $F$ (not necessarily finite)
such that $\overline{E}/\overline{F}$ is separable, and again put
$G = \Gal(E/F)$.
(In particular, if $\overline{F}$ is perfect, we may take $E = F^{\sep}$.)
Put $G^i = \varprojlim \Gal(E'/F)^i$,
where $E'$ runs over all finite Galois subextensions of $E/F$. By
\cite[Proposition~IV.14]{serre}, $\Gal(E'/F)^i$ is the image of $G^i$
in $\Gal(E'/F)$. Again, we call the resulting
filtration the \emph{upper numbering
filtration} of $G$; 
it satisfies the left continuity property $G^i = \cap_{j<i} G^j$
\cite[Remark 1, p.\ 75]{serre}. On the other hand, the upper numbering
filtration is not right continuous; define $G^{i+} = \cup_{j > i} G^j$,
which is a closed subgroup of $G^i$ (for the profinite topology) which
may be strictly smaller than $G^i$.
\end{defn}

\begin{defn}
With notation as in Definition~\ref{D:upper2},
we say that $i \geq 0$ is a \emph{break} of $E/F$ if
$G^i \neq G^{i+}$. (The term ``break'' is short for ``ramification break'';
the term ``jump'', for ``ramification jump'', is also used.)
If $E/F$ is finite, then there are finitely
many breaks; we refer to the largest one of them (or $0$ if there
are no breaks) as the \emph{highest break} of $E/F$, and denote it by
$b(E/F)$.
By Herbrand's theorem, if $E/F$ is the compositum of $E_1/F$ and $E_2/F$,
then
\[
b(E/F) = \max\{b(E_1/F), b(E_2/F)\}.
\]
The highest break is always rational (because it is the image of
an integer under $\phi_{E/F}$, which is a piecewise linear function
with rational slopes, breaks, and $y$-intercept) but is not necessarily
an integer. However, if $E/F$ is abelian, then the Hasse-Arf theorem asserts
that $b(E/F)$ is an integer \cite[Theorem~V.1]{serre}.
\end{defn}

\begin{defn} \label{D:not Galois}
Let $E/F$ be a finite but not necessarily Galois extension,
and let $E'/F$ be a finite
Galois extension containing $E$. From the transitivity of the
$\psi$ and $\phi$ functions (Remark~\ref{R:trans}), 
it follows (as in
\cite[Remark 2, p.~75]{serre}) that the function
$\phi_{E'/F} \circ \psi_{E'/E}$ depends only on $E$ and $F$ and 
not on $E'$. We call this function $\phi_{E/F}$; it has the property
that for $m > b(E/F)$ and $E'/F$ finite Galois containing $E$,
$b(E'/E) = m$ if and only if $b(E'/F) = \phi_{E/F}(m)$.
\end{defn}

\begin{remark}
If $E$ is a finite separable but not necessarily Galois extension of $F$,
we can define the highest break of $E$ to be the highest break of
its Galois closure.
If $E$ is not a field but only a finite \'etale $K$-algebra, then it
decomposes as a product of finite field extensions of $K$, 
and we can define the highest break $b(E/F)$
of $E$ to be the maximum of the highest breaks of any component of $E$.
With this convention, one has the rules
\begin{align*}
b(E_1 \oplus E_2/F) &= \max\{b(E_1/F), b(E_2/F)\} \\
b(E_1 \otimes E_2/F) &= \max\{b(E_1/F), b(E_2/F)\}.
\end{align*}
\end{remark}

\begin{example} \label{E:as}
For $k$ a field of characteristic $p>0$,
$F = k((t))$, and $E = k((t))[z]/(z^p - z - P(t^{-1}))$,
 where
$P$ is a polynomial over $k$ whose degree $d$ is not divisible by $p$,
a straightforward calculation \cite[Exercise~IV.2.5]{serre}
shows that $E/F$ has exactly one break, which is equal to $d$.
It follows that
\[
\phi_{E/F}(m) = \begin{cases} m & m \leq d \\
d + \frac{m-d}{p} & m > d.
                \end{cases}
\]
\end{example}

\begin{remark}
Note that for any finite separable extension $E$ of $F$,
$\phi_{E/F}$ is monotone: this follows for $E/F$
Galois by Definition~\ref{D:upper}, and for $E/F$ general by
Definition~\ref{D:not Galois}.
\end{remark}

\subsection{Unramified and tame extensions}

We next recall some more facts about extensions of local fields from
\cite{serre}, and extend a few definitions to the case of an
inseparable residue field extension. 
We retain all definitions and notations from the
previous section; we also continue to assume that $E/F$ is
an extension of complete discretely valued fields.

See
\cite[III.5]{serre} for all results implicit in the following
definitions.
\begin{defn} \label{D:unramified}
If $E/F$ is finite, 
we say $E/F$ is \emph{unramified} if $\gothm_E = \gothm_F \gotho_E$
and $\overline{E}/\overline{F}$ is separable.
Every subextension of an unramified extension is unramified,
so we may extend the definition to $E/F$ infinite
by saying that $E/F$ is unramified if every finite
subextension of $E/F$ is unramified in the previous sense.
\end{defn}

\begin{defn}
The compositum of unramified extensions is again unramified,
so any separable algebraic extension $E/F$
admits a \emph{maximal unramified subextension} $U$;
we say $E/F$ is \emph{totally ramified} if $U =F$.
By Hensel's lemma, $\overline{U}$ is the maximal separable subextension
of $\overline{E}/\overline{F}$. In particular, an unramified extension
is uniquely determined by its (separable) residue field extension;
if $F = k((t))$, this means that any finite unramified extension has the form
$k'((t))$ for some finite separable extension $k'/k$.
By \cite[Proposition~IV.2]{serre}, if $E/F$ is Galois with group $G$
and $\overline{E}/\overline{F}$
is separable, then the maximal unramified subextension of $E/F$ is 
the fixed field of $G_0 = G^0$.
\end{defn}

Oddly, the following quite standard definition does not occur in \cite{serre}.
\begin{defn}
If $E/F$ is finite Galois, we say $E/F$ is \emph{tamely ramified}
(or simply \emph{tame}) if $\Gal(E/U)$ has order coprime to $p$,
where $U$ is the maximal unramified subextension of $E/F$. 
Any subextension of a tame extension is tame, so we may extend
the definition to $E/F$ infinite by saying that $E/F$ is tame 
if each of its finite subextensions is tame. Also,
the compositum of tame extensions is tame, so any
Galois algebraic extension $E/F$ has a maximal tame subextension $T$.
If $\overline{E}/\overline{F}$ is separable, then
by the properties of the ramification filtration
stated in Definition~\ref{D:lower}, $T$ is the fixed field of $G^{0+}$.
(Note that Remark~\ref{R:tame} below implies that
if $E/F$ is itself tame, then $\overline{E}/\overline{F}$
is always separable.)
We say $E/F$ is \emph{totally wildly ramified} if $T = F$.
\end{defn}

\begin{remark} \label{R:tame}
If $E/F$ is finite Galois, totally ramified, and tame of degree $d$,
and moreover $F$ contains a primitive $d$-th root of unity $\zeta_d$, then
by Kummer theory \cite[X.3]{serre},
$T = F(\pi^{1/d})$ for some generator $\pi$ of $\gothm_F$. 
It follows in this case that $\overline{E} = \overline{F}$. If on the other
hand $\zeta_d \notin F$, then
$E(\zeta_d)$ and $F(\zeta_d)$ have the same residue field,
so $\overline{E}/\overline{F}$ is at least separable; since
$E/F$ was assumed to be totally ramified, we must have $\overline{E} = 
\overline{F}$. Finally, if $E/F$ is finite Galois and tame, and
$U$ is the maximal unramified subextension of $E/F$, then the previous
argument shows that $\overline{E} = \overline{U}$, so
$\overline{E}/\overline{F}$ is separable.
As noted in 
Definition~\ref{D:lower}, it then follows that 
$\Gal(E/U) = G_1/G_0$ is cyclic.
\end{remark}

\begin{defn}
If $E/F$ is finite Galois with maximal unramified subextension $U$ and
maximal tame subextension $U$, we define the \emph{tame degree} of $E/F$
to be the degree of $T$ over $U$, and the
\emph{wild degree} of $E/F$ to be the degree of $E$ over $T$.
Then the tame degree is coprime to $p$, and the wild degree is
a power of $p$.
\end{defn}

\subsection{Break decompositions}

We now recall some terminology regarding representations of the
absolute Galois group of a local field, following
\cite[Chapter~1]{katz-gauss}.

\begin{hypo} \label{H:breaks}
Throughout this section,
let $F$ be a complete discretely valued field whose residue field
$\overline{F}$ is perfect of characteristic $p>0$, 
put $G = \Gal(F^{\sep}/F)$, and put $P = G^{0+}$.
Then $P$ is the $p$-Sylow subgroup
of $G$ (in the sense of profinite groups).
\end{hypo}

\begin{convention}
When a group $G$ acts on a set $M$, let $M^G$ denote the fixed
set of $M$ under $G$.
\end{convention}

\begin{defn} \label{D:breaks}
Let $M$ be a $\ZZ[1/p]$-module on which $P$ acts via a finite
discrete quotient. Then by \cite[Proposition~1.1]{katz-gauss},
there is a unique direct sum decomposition $M = \oplus_{i \geq 0} M(i)$
of $M$ into $P$-stable submodules such that
\begin{align*}
M(0) &= M^{P} \\
M(i)^{G^i} &= 0 \qquad (i > 0) \\
M(i)^{G^j} &= M(i) \qquad (j>i).
\end{align*}
This decomposition is called the \emph{break decomposition} of $M$;
the associated descending filtration $M_i = \oplus_{j \geq i} M(j)$
is called the \emph{break filtration} of $M$.
There are finitely many $i \geq 0$ for which $M(i) \neq 0$,
and they are all rational numbers;
they are called the \emph{breaks} of $M$. If $M$ is nonzero, there must
be at least one break; the largest one is called the \emph{highest break}
of $M$, and is denoted $b_F(M)$.
\end{defn}

\begin{remark}
If $E$ is the fixed field of the kernel of the action of $P$ on $M$, then
the highest break of $M$ coincides with the highest break of $E/F$.
\end{remark}

Often instead of the full break decomposition, we consider some
of its numerical invariants. 
\begin{defn} \label{D:Hasse-Arf}
With notation as in Definition~\ref{D:breaks}, suppose that $M$ is
a free module over some ring $A$ and that $P$ acts $A$-linearly. Then
for each $i \geq 0$, $M(i)$ is projective of some finite rank; that
rank is called the \emph{multiplicity} of $i$ (as a break of $M$).
Define the \emph{Hasse-Arf polygon} of $M$, denoted $\calP(M)$,
as the polygon with left endpoint $(0,0)$ consisting of $n$ segments of 
horizontal width 1, the $i$-th of which has slope equal to the $i$-th
smallest break of $\rho$, counting multiplicities. 
\end{defn}

The strong form of the Hasse-Arf theorem
\cite[Proposition~1.9]{katz-gauss} yields the following integrality property
of the Hasse-Arf polygon.
\begin{prop} \label{P:Hasse-Arf}
With notation as in Definition~\ref{D:Hasse-Arf}, the Hasse-Arf
polygon has integer vertices. In particular, the breaks of any abelian
extension of $F$ are all integers.
\end{prop}

\begin{remark}
While the Hasse-Arf polygon looks like a Newton polygon of the sort
one associates to a polynomial over a local field, its formalism
is quite different.
For instance, the highest break of the tensor product of two modules
is at most the maximum of the highest breaks of the tensorands, whereas
the highest slope of the tensor product of two polynomials (that is,
the polynomial whose roots are the products of roots, one from
each tensorand) is the sum of the highest slopes of the tensorands.
For a thorough development of the formalism of Hasse-Arf polygons
(in which it is shown that any class of filtrations
which ``look enough like'' ramification filtrations actually
\emph{are} ramification filtrations), see \cite{andre}.
\end{remark}

%\begin{defn}
%Suppose $M, M'$ are two modules with $P$-actions as in
%Definition~\ref{D:Hasse-Arf}. We write $\calP(M) \leq \calP(M')$
%to mean that $\calP(M)$ lies
%entirely on or below $\calP(M')$, or equivalently, that
%$\rank(M_i) \leq \rank(M'_i)$ for all $i \geq 0$.
%\end{defn}

Using the Hasse-Arf theorem, we can obtain the following finiteness
result about representations of $P$. First, recall Jordan's theorem
on finite linear groups.
\begin{prop} \label{P:Jordan}
For any positive integer $n$, there exists an integer $f(n)$ such that
for any field $K$ of characteristic zero, 
any finite subgroup $G$ of $\GL_n(K)$ contains a commutative normal subgroup
$H$ of index at most $f(n)$.
\end{prop}
\begin{proof}
Any such $G$ can be embedded into $\GL_n(\QQ^{\alg})$, and hence into
$\GL_n(\CC)$. For the result in this case, see \cite{isaacs}.
\end{proof}

\begin{prop} \label{P:bound image}
Given a residual characteristic
$p$, a positive integer $n$ and a nonnegative real number $\ell$, there
exists an integer $N = N(p,n,\ell)$ such that every representation of $P$
of dimension $n$ over a field of characteristic zero, with
finite discrete image and
highest break $\leq \ell$, has image of order at most $N$.
\end{prop}
\begin{proof}
There is no loss of generality in working with representations
over $\CC$ (or even over the algebraic closure of $\QQ$). 
Let $\rho: P \to \GL_n(\CC)$ be a representation with image $G
\subset \GL_n(\CC)$, and let $E$ be the fixed field of the kernel
of $\rho$.
Suppose first that $G$ is abelian; then
by Proposition~\ref{P:Hasse-Arf}, all of the breaks of $E/F$
are integers, so the number of them is bounded by $\ell$.
Moreover, an elementary abelian $p$-subgroup of $\GL_n(\CC)$ can
have at most $p^n$ elements since the matrices in such a subgroup
must be simultaneously diagonalizable. Thus $G$ has order at
most $p^{n\ell}$ in this case.

Now suppose that $G$ is arbitrary.
Apply Jordan's theorem (Proposition~\ref{P:Jordan})
to choose a commutative normal subgroup $H$ of $G$ of index bounded
as a function of $n$,
and let $E$ be the fixed field of $H$; then the highest break
of the restriction of $\rho$ to $P \cap \Gal(E^{\sep}/E)$
is bounded by a function of $p,n,\ell$. This restriction is abelian,
so as before, the order of $H$ is bounded by a function of $p,n,\ell$.
This yields the desired result.
\end{proof}

\section{Rigid annuli}
\label{sec:rigid}

We now introduce the rigid analytic spaces we will be working with,
which are certain one-dimensional annuli. In particular, we verify that
over a spherically complete coefficient field, every coherent
locally free sheaf on a one-dimensional rigid annulus is freely
generated by global sections (Theorem~\ref{T:coherent}).
This provides a bridge between our setup and the existing
literature on $p$-adic differential equations, which is mostly conducted
in ring-theoretic terms.
We also produce a special class of finite \'etale covers
of rigid annuli corresponding to finite \'etale extensions of
$k((t))$; these will be used to discuss the monodromy of
$p$-adic differential equations.

We will freely use the language of rigid analytic geometry
using the original foundations of Tate et al; see \cite{fresnel} (particularly
Chapter~2) for an introduction. 
If one prefers the Berkovich foundations, as in
\cite{berkovich-ihes} (see also \cite{berkovich-icm} for an overview), 
one should in principle have no trouble converting
the discussion into those terms, since the rigid spaces under consideration
are quasi-separated and admit affinoid coverings of finite type.
However, one will probably encounter some subtleties; see for instance
Remark~\ref{R:continuous}.

\subsection{Notations}
\label{subsec:notations}

Before proceeding, we set a few notational conventions.

\begin{convention} \label{conv:field}
Let $K$ be a field complete with respect to a nonarchimedean absolute
value $|\cdot|: K^* \to \RR^+$.
Let $\gotho_K$ be the subring of $x \in K$ with $|x| \leq 1$,
let $\gothm_K$ be the ideal of $x \in \gotho_K$ with $|x| < 1$,
and let $k$ denote the residue field $\gotho_K/\gothm_K$.
Let $\Gamma^*$ denote the divisible closure of the image of $|\cdot|$. 
\end{convention}

\begin{convention}
On any rigid analytic space over $K$, let $\calO$ denote the structure
sheaf, and let $\gotho$ denote the subsheaf of the structure sheaf
consisting of functions bounded in absolute value by 1 everywhere
(i.e., the ``integral subsheaf'' of $\calO$).
\end{convention}

\begin{convention}
In case the notation for an object includes an explicit mention of the
coefficient field $K$, we will routinely suppress that $K$ from the
notation when the choice of $K$ is to be understood (as it almost always
will be, except when we need to compare distinct choices).
For example, in Definition~\ref{D:Robba}, we typically abbreviate
$\calR_{I,K}$ to $\calR_I$.
\end{convention}

\begin{convention}
When we write the absolute value of a matrix, we mean the maximum
of the absolute values of its entries, rather than any sort of
operator norm. (Contrast this convention with Definition~\ref{D:operator
norm}.)
Also, we let $I_n$ denote the $n \times n$ identity matrix over any ring.
\end{convention}

\subsection{Rigid annuli}

\begin{defn} \label{D:Robba}
We say a subinterval $I$ of $[0,1)$ is \emph{aligned} if any nonzero
endpoint at which it is closed is contained in $\Gamma^*$.
For $I$ aligned, we define the annulus $A(I) = A_K(I)$ as the admissible
open subspace of the rigid affine $t$-line given by
\[
\{t \in \AAA^1_K: |t| \in I\}.
\]
If $I$ is given with explicit endpoints, the enclosing
parentheses are omitted, so that we write for instance
$A[0,1)$ instead of $A([0,1))$.
Let $\calR_I = \calR_{I,K}$
denote the ring $\Gamma(\calO,A(I))$ of rigid analytic functions
on $A(I)$.
The elements of $\calR_I$ can be described as formal
Laurent series $\sum_{i \in \ZZ} c_i t^i$ with each $c_i \in K$; 
for $r \in I \cap \Gamma^*$, the spectral seminorm on the subspace $|t| = r$
of $A(I)$, restricted to $\calR_I$, is equal to the
norm $|\cdot|_r$ given by the formula
\[
\left|\sum c_i t^i \right|_r = \sup_i \{|c_i| r^i\}.
\]
\end{defn}

One has analogues of the maximum modulus principle and the Hadamard 
three circles theorem for $|\cdot|_{r}$.
\begin{lemma} \label{L:Hadamard}
\begin{enumerate}
\item[(a)]
For $x \in \calR_{[0,b]}$ and $r \in [0,b]$,
$|x|_{r} \leq |x|_{b}$.
\item[(b)]
For $x \in \calR_I$, $a,b \in I$, and $c \in [0,1]$,
put $r = a^{c} b^{1-c}$; then
$|x|_{r} \leq |x|_{a}^c |x|_{b}^{1-c}$.
\end{enumerate}
\end{lemma}
\begin{proof}
\begin{enumerate}
\item[(a)]
If $x = \sum c_j t^j \in \calR_{[0,b]}$,
then $c_j = 0$ for $j<0$.
Hence if $r \in [0,b]$, then
\[
|c_j| r^j \leq |c_j| b^j;
\]
taking suprema yields $|x|_{r} \leq |x|_{b}$.
\item[(b)]
Note that the desired inequality holds with equality if
$x = c_j t^j$.
For a general $x = \sum c_j t^j$, we then have
\begin{align*}
|x|_{r} &= \sup_j \{ |c_j t^j|_{r} \} \\
&\leq \sup_j \{ |c_j t^j|_{a}^c |c_j t^j|_{b}^{1-c} \} \\
&\leq \sup_j \{ |c_j t^j|_{a} \}^c \sup_j \{|c_j t^j|_{b} \}^{1-c} \\
&= |x|_{a}^c |x|_{b}^{1-c},
\end{align*}
as desired.
\end{enumerate}
\end{proof}

A result of Lazard \cite[Proposition~4]{lazard} yields the following.
\begin{prop} \label{P:poly generator}
Let $I$ be a closed aligned subinterval of $[0,1)$. Then every
ideal of $\calR_I$ is generated by an element
of $K[t]$. In particular, $\calR_I$ is a principal
ideal domain.
\end{prop}

\begin{remark}
The alignedness restriction can be dumped if one allows affinoid spaces
to be closed analytic subspaces of polydiscs of arbitrary radii, not 
just radius 1. This permission is made in Berkovich's foundations
of rigid geometry.
\end{remark}

\subsection{A matrix approximation lemma}

We will need a matrix approximation lemma in the spirit of
\cite[Lemma~6.2]{me-local}.
We start with an analogue of \cite[Lemma~6.1]{me-local}.
\begin{lemma} \label{L:add lambda}
Put $S = K[t]$ or $K[t,t^{-1}]$, and fix $x,y \in S$.
Then there exists $c>0$
such that for any $\lambda \in K$ with $|\lambda| < c$, 
$x$ and $y+\lambda$ generate the unit ideal in $S$.
\end{lemma}
\begin{proof}
For each $\lambda \in K$, the ideal generated by $x$ and $y + \lambda$
in $S$ is generated by some monic polynomial $e_\lambda$.
Note that this limits $e_\lambda$ to a finite set, namely the
monic factors of $x$. Each value of $e_\lambda$ not equal to 1 can
only occur for one value of $\lambda$: if $e_\lambda = e_{\lambda'}$,
then $e_\lambda$ divides $(y+\lambda) - (y+\lambda') = 
\lambda - \lambda' \in K$, contradiction.
In particular, for $c>0$ sufficiently small,
$e_\lambda = 1$ for all $\lambda \in K$ with $|\lambda| < c$,
as desired.
\end{proof}

The following is analogous to \cite[Lemma~6.2]{me-local}, but with
some slight simplifications because we are taking $I$ to be closed,
and because we are only working with power series rings instead
of the more general ``analytic rings'' of \cite{me-local}.
\begin{lemma} \label{L:approximate matrix}
Let $I$ be a closed aligned subinterval of $[0,1)$ which
does (resp.\ does not) contain $0$,
and let $M$ be an invertible $n \times n$ matrix over $\calR_I$.
Then there exists an invertible $n \times n$ matrix $U$ over
$K[t]$ (resp.\ over $K[t,t^{-1}]$) such that $|MU-I_n|_r < 1$ for $r \in I$.
Moreover, if $|\det(M)-1|_r < 1$, we can ensure that $\det(U) = 1$.
\end{lemma}
\begin{proof}
We proceed by induction on $n$, the case $n=1$ being vacuously true.
By multiplying some row of $M$ by the inverse of
$\det(M)$, we may reduce to the case where
$\det(M) = 1$.
Let $C_i$ denote the cofactor of $M_{ni}$ in $M$, so that
$\det(M) = \sum_{i=1}^n C_i M_{ni}$, and in fact
$C_i = (M^{-1})_{in} \det(M)$.
Thus $C_1, \dots, C_n$
generate the unit ideal in $\calR_I$, so we can find
$\alpha_1, \dots, \alpha_n \in \calR_I$ such that
$\sum_{i=1}^n \alpha_i C_i = 1$.

For brevity, write $S = K[t]$ if $0 \in I$ and $S = K[t,t^{-1}]$ if 
$0 \notin I$. Choose $\beta_1, \dots, \beta_{n-1}, \beta'_n \in S$
such that for $r \in I$,
\[
|\beta_i - \alpha_i|_r < \min_j \{|C_j|^{-1}_r\} \quad (i=1, \dots, n-1),
\qquad
|\beta'_n - \alpha_n|_r < \min_j \{|C_j|^{-1}_r\}.
\]
By Lemma~\ref{L:add lambda}, for $\lambda \in K$ of sufficiently
small absolute value, $\beta_n = \beta'_n + \lambda$ has the properties
that $|\beta_n - \alpha_n|_r < \min_j\{|C_j|_r\}$ for $r \in I$,
and that $\beta_1, \dots, \beta_n$ generate the unit ideal in
$S$. Since $S$ is a principal ideal domain, we can find a matrix
$A$ over $S$ of determinant 1 such that
$A_{ni} = \beta_i$ for $i=1, \dots, n$.
Put $M' = MA^{-1}$, and let $C'_n$ be the cofactor of
$M'_{nn}$ in $M'$. Then
\begin{align*}
C'_n &= (M')^{-1}_{nn} \det(M') \\
&= (AM^{-1})_{nn} \det(M) \\
&= \sum_{i=1}^n A_{ni} M^{-1}_{in} \det(M) \\
&= \sum_{i=1}^n \beta_i C_i,
\end{align*}
so that
\[
C'_n = 1 + \sum_{i=1}^n (\beta_i - \alpha_i) C_i
\]
and so $|C'_n - 1|_r < 1$ for $r \in I$. In particular, $C'_n$
is a unit in $\calR_I$.

Apply the induction hypothesis to the upper left $(n-1)\times (n-1)$
submatrix of $M'$, and extend the resulting matrix $V$ to
an $n \times n$ matrix by setting $V_{ni} = V_{in} = 0$
for $i=1, \dots, n-1$ and $V_{nn} = 1$. Then we have
$\det(M'V) = 1$ and 
\[
|(M'V - I_n)_{ij}|_r < 1 \qquad (i=1,\dots,n-1; \, j=1,\dots, n-1; \, r \in I).
\]

We now perform an ``approximate Gaussian elimination'' over $\calR_I$ 
to transform
$M'V$ into a new matrix $N$ with $|N-I_n|_r < 1$ for $r \in I$.
First, define a sequence of matrices $\{X^{(h)}\}_{h=0}^\infty$ by
$X^{(0)} = M'V$ and 
\[
X_{ij}^{(h+1)} = \begin{cases} X_{ij}^{(h)} & i<n \\
X_{nj}^{(h)} - \sum_{m=1}^{n-1} X_{nm}^{(h)} X_{mj}^{(h)} & i=n;
                 \end{cases}
\]
note that $X^{(h+1)}$ is obtained from $X^{(h)}$ by subtracting
$X_{nm}^{(h)}$ times the $m$-th row from the $n$-th for 
$m=1,\dots,n-1$ in succession. At each step, for each $r \in I$,
$\max_{1 \leq j \leq n-1} \{|X_{nj}^{(h)}|_r\}$
gets multiplied by a factor no larger than $\max_{1 \leq i,j \leq n-1}
\{|(M'V - I_n)_{ij}|_r\}$; since $I$ is closed and $|\cdot|_r$ is 
a continuous function of $r$, these factors are
bounded strictly below 1. 
Thus for $h$ sufficiently large, we have
\[
|X_{nj}^{(h)}|_r < \min\{1, \min_{1 \leq i \leq n-1} \{|X_{in}^{(h)}|^{-1}
\}\} \qquad (r \in I;\, j = 1 ,\dots, n-1).
\]
Pick such an $h$ and set $X = X^{(h)}$; note that $\det(X) = \det(M'V) = 1$. 
Then for $r \in I$,
\begin{align*}
|(X-I_n)_{ij}|_r &< 1 \qquad (i=1,\dots, n; \,j =1, \dots, n-1) \\
|X_{in} X_{nj}|_r &<1 \qquad (i=1,\dots, n-1; \,j =1, \dots, n-1) 
\end{align*}
and hence also $|X_{nn} - 1|_r < 1$.

Next, define a sequence of matrices $\{W^{(h)}\}_{h=0}^\infty$ by
setting $W^{(0)} = X$ and
\[
W_{ij}^{(h+1)} = \begin{cases} W_{ij}^{(h)} - W_{in}^{(h)} W_{nj}^{(h)} &
i<n \\
W_{ij}^{(h)} & i=n; \end{cases}
\]
note that $W^{(h+1)}$ is obtained from $W^{(h)}$ by subtracting
$W_{in}^{(h)}$ times the $n$-th row from the $i$-th row for
$i=1, \dots, n-1$. At each step, for $r \in I$,
$|W_{in}^{(h)}|_r$ gets multiplied by a factor no larger than
$|X_{nn}^{(h)}-1|_r$; again, these factors are bounded strictly below 1
because $I$ is closed. 
Thus for $h$ sufficiently large,
\[
|W^{(h)}_{in}|_r < 1 \qquad (r \in I;\, 1 \leq i \leq n-1).
\]
Pick such an $h$ and set $W = W_h$; then $|W-I_n|_r < 1$ for $r \in I$.
(Note that the inequality $|X_{in} X_{nj}|_r < 1$ for
$i=1, \dots, n-1$ and $j=1, \dots, n-1$ ensures that
the second set of row operations does not disturb the fact that
$|W^{(h)}_{ij}|_r < 1$ for $i=1,\dots, n-1$ and $j=1,\dots, n-1$.)

To conclude, note that by construction, $(M'V)^{-1} W$ is a product
of elementary matrices over $\calR_I$, each consisting of the diagonal
matrix plus one off-diagonal entry. By suitably approximating the
off-diagonal entry of each
matrix in the product by an element of $S$, we get an invertible
matrix $X$ over $S$ such that $|M'VX - I_n|_r < 1$ for $r \in I$.
We may thus take $U = A^{-1}VX$ to obtain the desired result.
\end{proof}
\begin{remark}
It is tempting to believe that one can improve the conclusion of
Lemma~\ref{L:approximate matrix} to yield that for any given $c>0$,
we can choose $U$ so that
$|MU-I_n|_r < c$ for $r \in I$, as in \cite[Lemma~6.3]{me-local}.
However, this is only possible if we assume
$|\det(M) - 1|_r < c$ for $r \in I$, or else the base case $n=1$
fails. Indeed, this failure is inevitable, since Theorem~\ref{T:coherent}
does not hold for arbitrary $K$; see Remark~\ref{R:not spherical}.
\end{remark}

\subsection{Locally free sheaves on rigid annuli}

Recall that a theorem of Kiehl \cite[Theorem~4.5.2]{fresnel} asserts
that a coherent sheaf on an affinoid space is generated by finitely
many global sections. This is typically not true on a nonaffinoid
space, such an as open rigid annulus, but for coherent locally
free sheaves, it turns out we can salvage something. First,
however, we must restrict the field of coefficients.

\begin{defn}
The field $K$, which we are supposing to be complete for
a nonarchimedean absolute value, is said to be \emph{spherically
complete}, or \emph{maximally complete},
if every decreasing sequence of closed balls
has nonempty intersection. For instance, every discretely valued
field is spherically complete, but the field $\CC_p$, the completed
algebraic closure of $\QQ_p$, is not spherically complete.
\end{defn}

\begin{prop} \label{P:spherical}
Suppose that $K$ is spherically complete.
Fix a sequence $(r_1, r_2, \dots)$ of positive real numbers, and
let $c_0$ be the set of sequences $x = (x_1, x_2, \dots)$ over $K$
with $|x_i| \leq r_i$ for all $i$.
For $i=0,1,\dots$, let $f_i$ be an affine functional on $c_0$, that is,
\[
f_i(x_1, x_2, \dots) = a_{i,0} + \sum_{j=1}^\infty a_{i,j} x_j
\]
for some sequence $a_{i,0},a_{i,1},\dots$ converging to $0$ in $K$.
Let $S_i$ be the subset of $x \in c_0$ on which $|f_h(x)| < 1$ for
$h=0, \dots, i$. If $S_i \neq \emptyset$ for each $i$, then
$\cap_i S_i \neq \emptyset$.
\end{prop}
\begin{proof}
Suppose that $y_1, y_2, \dots, y_l \in K$ have been chosen so that 
$|y_j| \leq r_j$ for $j=1, \dots, l$, and for each $i$, 
the set $S_{i,l}$ of $x \in c_0$ with $x_1 = y_1, \dots, x_l = y_l$ and
$|f_h(x)| < 1$ for $h=0, \dots, i$ is nonempty.
Let $T_{i,l}$ be the set of possible values of $x_{l+1}$ for a sequence
$x \in S_{i,l}$. Then $T_{i,l}$ is necessarily an open ball, and the sequence
$T_{0,l}, T_{1,l}, \dots$ is decreasing,
so has a nonempty intersection since $K$ is spherically
complete. We may then pick any $y_{l+1}$ in that intersection to continue the
construction.

The hypothesis in the previous paragraph holds vacuously with $l=0$.
We may thus construct $y_1,y_2, \dots$ as above, and the resulting sequence
belongs to $\cap_i S_i$.
\end{proof}

\begin{theorem} \label{T:coherent}
Let $I$ be an aligned subinterval of $[0,1)$ and let
$\calE$ be a coherent locally free sheaf of rank $n$ on $A(I)$. Then
there exist sections $\bv_1, \dots, \bv_n \in \Gamma(\calE, A(I))$
which freely generate $\calE$.
\end{theorem}
\begin{proof}
Let $J_1 \subseteq J_2 \subseteq \cdots$ be a weakly increasing
sequence of aligned closed subintervals of $I$ whose union is all of
$I$, with the property that if $0 \in I$, then $0 \in J_i$ for all $i$.
Put $\calR_i = \calR_{J_i}$ and $E_i = \Gamma(\calE,
A(J_i))$; by Proposition~\ref{P:poly generator}, each $E_i$
is free of rank $n$ over $\calR_i$.

Choose a basis $\bv_{1,1}, \dots, \bv_{1,n}$ of $E_1$.
Given a basis $\bv_{i,1}, \dots, \bv_{i,n}$ of $E_i$, we
choose a basis $\bv_{i+1,1}, \dots, \bv_{i+1,n}$ of $E_{i+1}$ as 
follows. Pick any basis $\be_1, \dots, \be_n$ of $E_{i+1}$, and define
an invertible $n \times n$ matrix $M_i$ over $\calR_i$ by writing
$\be_l = \sum_{j} (M_i)_{jl} \bv_{i,j}$. Apply Lemma~\ref{L:approximate
matrix} to produce an invertible $n \times n$ matrix $U$ over $S$,
where $S = K[t]$ if $0 \in I$ and $S = K[t,t^{-1}]$ if $0 \notin I$,
such that $|M_i U-I_n|_r < 1$ for $r \in J_{i}$. 
Put $V_i = M_i U$, and define the
basis $\bv_{i+1,1}, \dots, \bv_{i+1,n}$ of $E_{i+1}$ by
$\bv_{i+1,l} = \sum_j (V_i)_{jl} \bv_{i,j}$.

Now suppose $\be_1, \dots, \be_n$ is a basis of $E_1$. Define
the invertible $n \times n$ matrix $W$ over $\calR_1$ by
\[
\be_l = \sum_j W_{jl} \bv_{1,j}.
\]
Put $W_i = (V_1 \cdots V_{i-1})^{-1} W$; then
\[
\be_l = \sum_j (W_i)_{jl} \bv_{i,j}.
\]  
Hence $\be_1, \dots, \be_n$ forms a basis of $E_i$ if and only
if $W_i$ is an invertible matrix over $\calR_i$.

Let $B_i$ denote the set of $n \times n$ matrices $W$ over $\calR_1$ with
\[
|(V_1 \cdots V_{i-1})^{-1} W - I_n|_r < 1 \qquad (r \in J_{i});
\]
since $|V_i-I_n|_r < 1$ whenever $r \in J_i$, we have
$B_{i+1} \subseteq B_i$.
We may apply Proposition~\ref{P:spherical} by fixing some $s \in J_1$
and identifying
$c_0$ with the set of $n \times n$ matrices $W$, with
$W_{ij} = \sum_l W_{ij,l} t^l$, satisfying $|W_{ij,l} t^l|_s \leq 1$ for all 
$i,j,l$; by doing so, we see that $\cap_i B_i \neq \emptyset$.
If $W \in \cap_i B_i$, we may put
$\be_l = \sum_j W_{jl} \bv_{1,j}$ to obtain a basis
$\be_1, \dots, \be_n$ of $E_1$ that extends to a basis of $E_i$ for
each $i$. Thus the $\be_i$ are global sections of $\calE$ which
freely generate $\calE$, as desired.
\end{proof}

\begin{remark} \label{R:not spherical}
Theorem~\ref{T:coherent} is false for \emph{any} field $K$ which is not
spherically complete, even for line bundles. 
More precisely, Lazard \cite[Th\'eor\`eme~2]{lazard} showed that
on an open disc over $K$, every coherent locally free sheaf of rank 1
is generated by a global section (and hence trivial) if and only if
$K$ is spherically complete. 
Thus Theorem~\ref{T:coherent} may be viewed
as a higher-rank generalization of Lazard's result. (For $I = [0,1)$,
i.e., for locally free sheaves on an open disc,
this generalization was already given by Gruson
\cite[Proposition~2]{gruson}.) It may also be 
viewed as an explication of a special case of the comment made
in the introduction of \cite{vdp}, to the effect that the proof given
there that the sheaf $\calO^*$ has no higher cohomology on any rational
subset of the projective line can be carried over to establish
the vanishing of $H^1(\GL_n)$.
\end{remark}

\subsection{Robba rings and formally \'etale covers}

\begin{defn}
Let $\calR = \calR_K$ denote the direct limit of the rings $\calR_{(a,1)}$ over
all $a \in (0,1)$; the ring $\calR$ is called the
\emph{Robba ring} over $K$. The elements of $\calR$ can be viewed as 
formal Laurent series
$\sum_i c_i t^i$, with $c_i \in K$, which converge on some unspecified
open annulus with outer radius 1. Let $\calR^{\inte} = \calR_K^{\inte}$ denote
the subring of $\calR$ consisting of series with $c_i \in \gotho_K$
for all $i \in \ZZ$; let $\gothm_{\calR} = \gothm_{\calR_K}$ 
denote the ideal of
$\calR^{\inte}$ consisting of series with $c_i \in \gothm_K$
for all $i \in \ZZ$.
\end{defn}

\begin{lemma}
The ring $\calR^{\inte}$ is a local ring with maximal ideal
$\gothm_{\calR}$ and residue field $k((t))$.
\end{lemma}
\begin{proof}
For any $x \in \gothm_{\calR}$,
we have $|x|_a < 1$ for $a$ sufficiently close to 1, so
$1+x$ is invertible in $\calR^{\inte}$. Hence 
$\gothm_{\calR}$ is contained in the Jacobson radical of
$\calR^{\inte}$, but the Jacobson radical is the intersection of all
maximal ideals. Hence $\gothm_{\calR}$ is the unique maximal ideal
of $\calR^{\inte}$, proving the claim.
\end{proof}

\begin{remark}
If $K$ is discretely valued, then $\calR^{\inte}$ is a discrete
valuation ring with corresponding absolute
value $|\sum c_i t^i| = \sup_i \{|c_i|\}$. 
On the other hand, if $K$ is not discretely valued, then one can
construct $\sum c_i t^i \in \gothm_{\calR}$ with $\sup_i \{|c_i|\} = 1$,
so we cannot view $\calR^{\inte}$ as a valuation ring in this fashion.
See \cite[Remark~14]{christol} for an example of how
readily this discrepancy can crop up if one does not take pains to avoid it.
\end{remark}

\begin{defn}
A pair $(R,I)$, in which $R$ is a ring and $I$ is an ideal of $R$, is 
\emph{henselian} if each finite \'etale extension of $R/I$ lifts
uniquely to a finite \'etale extension of $R$.
\end{defn}
\begin{prop} \label{P:henselian}
The pair $(\calR^{\inte}, \gothm_{\calR})$ is henselian.
\end{prop}
\begin{proof}
By \cite[43.2]{nagata}, it suffices to show that for any monic polynomial
$P(x) = x^n + c_{n-1}
x^{n-1} + \cdots + c_0$ over $\calR^{\inte}$ with the property that
$c_{n-1}+1 \in \gothm_{\calR}$ and $c_i \in \gothm_{\calR}$ for 
$i=0, \dots, n-2$,
there exists a root $z \in \calR^{\inte}$ of $P(x)$ such that
$z-1 \in \gothm_{\calR}$. We construct $z$ using a Newton iteration as
follows.

Put $z_0 = 1$; given $z_i \in \calR^{\inte}$ such that
$z_i - 1 \in \gothm_{\calR}$, note that
\begin{align*}
P'(z_i) - 1 &= n z_i^{n-1} + (n-1) c_{n-1} z_i^{n-2} + \cdots + c_1\\
&\equiv n + (n-1) c_{n-1} \\
&\equiv 1 \pmod{\gothm_{\calR}},
\end{align*}
so $P'(z_1)$ is a unit in $\calR^{\inte}$. We may then put
\[
z_{i+1} = z_i - \frac{P(z_i)}{P'(z_i)};
\]
since $P(z_i) \equiv z_i^n + c_{n-1} z_i^{n-1} \equiv 0 \pmod{\gothm_{\calR}}$,
we have $z_{i+1} - 1 \in \gothm_{\calR}$ and the iteration continues.

Choose $a \in (0,1) \cap \Gamma^*$ such that for $r \in [a,1)$,
$|c_{n-1}+1|_r < 1$ and $|c_i|_r < 1$ for $i=0, \dots, n-2$.
Then by induction, for $r \in [a,1)$
we have $|z_i-1|_r < 1$. Moreover, for $r \in [a,1)$, we have
\begin{align*}
|P(z_{i+1})|_r &= \left|\sum_{j=2}^n \frac{P^{(j)}(z_i)}{j!}
(z_{i+1} - z_i)^j \right|_r \\
&\leq |z_{j+1}-z_j|_r^2,
\end{align*}
so $|P(z_i)|_r \to 0$ and hence $|z_{i+1}-z_i|_r \to 0$. It follows
that the $z_i$ converge to a limit $z$ satisfying $|z-i|_r < 1$ for
$r \in [a,1)$ and $P(z) = 0$, as desired.
\end{proof}

\begin{defn}
We say a finite cover $X$ of $A[a,1)$ is \emph{formally \'etale}
if it induces a finite \'etale extension of $\calR^{\inte}$;
we refer to the induced finite \'etale extension of 
$\calR^{\inte}/\gothm_{\calR}
\cong k((t))$ as the \emph{reduction} of $X$. 
By Proposition~\ref{P:henselian}, any two formally \'etale covers of $A[a,1)$
with isomorphic reductions
become themselves isomorphic 
over $A[b,1)$ for some $b \in [a,1) \cap \Gamma^*$.
\end{defn}

\begin{remark}
If the reduction of a formally \'etale cover $X \to A[a,1)$ 
induces a separable residue field extension of $k$, then $X$ itself
is isomorphic to an annulus over some finite extension of $K$ which
induces an \'etale extension of $\gotho_K$.
\end{remark}

\section{Monodromy of differential equations}
\label{sec:mono}

We next collect some facts about $p$-adic differential equations
on rigid annuli, specifically isolating those with quasi-unipotent
monodromy. Our treatment follows somewhat that of Matsuda
\cite{matsuda-katz}, though we simplify his presentation a bit by
making more systematic use of local duality (as in \cite{crewfin}).

\setcounter{equation}{0}
\begin{convention} \label{conv:mono}
Retain the notations introduced in Section~\ref{subsec:notations},
but now assume further that
$\charac(K) = 0$,
that $\charac(k) = p>0$, and that the absolute value $|\cdot|$ on $K$ is
normalized so that $|p| = p^{-1}$.
\end{convention}

\subsection{$\nabla$-modules on annuli: generalities}

To begin with, we define the category in which we will be working.
\begin{defn}
For $r \in (0,1) \cap \Gamma^*$, let $\calM_r = \calM_{r,K}$ 
denote the category of
$\nabla$-modules, i.e., coherent locally free sheaves $\calE$
equipped with a connection $\nabla: \calE \to
\calE \otimes \Omega^1$, on the annulus $A[r,1)$. Then there are natural
restriction functors $\calM_{r} \to \calM_{s}$ whenever $s \in [r,1) \cap 
\Gamma^*$;
using these, we may define the direct limit of the $\calM_{r}$, which
we denote by $\calM = \calM_K$.
\end{defn}

\begin{convention} \label{conv:shrinking}
Unless otherwise specified, any $\nabla$-module defined over
an annulus $A[r,1)$ will be interpreted as an object of
$\calM$; that is, a ``morphism'' between two such objects need
only be defined on some possibly smaller annulus $A[s,1)$,
and so on.
\end{convention}

\begin{remark}
On higher-dimensional spaces, one would ordinarily require the connection
$\nabla$ to be integrable, i.e., the composition of $\nabla$ with the
induced map $\calE \otimes \Omega^1 \to \calE \otimes \Omega^2$
should vanish. This is superfluous here because we are working on 
a one-dimensional space, so $\Omega^1$ is freely generated by $dt$ and
$\Omega^2$ vanishes.
\end{remark}

\begin{remark}
To specify a connection on a sheaf $\calE$,
it is enough to specify the action on $\calE$ of any
one differential operator of $A(I)$, e.g.,
$\frac{d}{dt}$ or $t \frac{d}{dt}$. Again, this is because
the underlying space is one-dimensional; on a $d$-dimensional space,
one must specify the actions of $d$ operators, and integrability
is equivalent to the fact that these actions commute.
\end{remark}

\begin{remark}
Note that given a $\nabla$-module $\calE$ on $A[r,1)$, any
$\nabla$-submodule $\calF$ is a locally free subsheaf; that is,
the quotient $\calE/\calF$ is torsion-free. This is a standard property
of $\nabla$-modules on any smooth rigid analytic space; see for instance
\cite[Proposition~2.2.3]{ber2}.
\end{remark}

We will occasionally need to enlarge the coefficient field $K$, as follows.
\begin{defn}
By an ``unramified (algebraic) extension of $K$'', we will mean an algebraic
extension $K'$ of $K$ such that the integral closure $\gotho_{K'}$
of $\gotho_K$ in $K'$ is \'etale over $\gotho_K$. This is inconsistent
with scheme-theoretic terminology, but is consistent with
the established terminology for local fields
(Definition~\ref{D:unramified}) in case $K$ is discretely valued.
Let $K^{\unr}$ denote the maximal unramified extension of $K$
(within $K^{\sep}$).
\end{defn}

%\begin{defn}
%We say that $\calE \in \calM$ is \emph{geometrically
%irreducible} if $\calE$ is irreducible as an element of
%$\calM_{K'}$ for any finite extension $K'$ of $K$.
%\end{defn}

\begin{defn}
We say that $\calE \in \calM$ is 
said to be \emph{overconvergent} if for any $\eta \in (0,1)$, there
exists $r \in (0,1)$ such that for any closed aligned subinterval
$I$ of $[r,1)$ and any $\bv \in \Gamma(\calE, A(I))$, the sequence
\[
\left\{ \frac{1}{j!} \frac{d^j}{dt^j} \bv \right\}_{j=0}^\infty
\]
is $\eta$-null; that is, in terms of some basis, the norms of the
coefficients of the $j$-th term of the series, multiplied by $\eta^j$,
converge to zero. It suffices to check this condition for some set
of intervals $I$ covering $[r,1)$, as then the corresponding subannuli
form an admissible cover of $A(I)$. 
Let $\calM^{\conv} = \calM^{\conv}_K$ be the subcategory of $\calM$ consisting
of overconvergent objects; this category is an abelian subcategory of
$\calM$.
\end{defn}
\begin{remark}
This notion of 
``overconvergent'' coincides with the notion of ``solvable at $1$'' in 
the terminology of \cite[4.1-1]{cm3}.
The property of overconvergence can be shown by a direct calculation
to be invariant under
automorphisms of $A[0,1)$.
\end{remark}

\begin{remark} \label{L:convergent pullback}
Let $f_1,f_2: A[a,1) \to A[b,1)$ be two morphisms which induce
the same map on $\calR^{\inte}/\gothm_{\calR}$. Then there is
a natural isomorphism between the pullback functors $f_1^*$ and
$f_2^*$ on $\calM^{\conv}$, given by the Taylor series
\[
1 \otimes \bv \mapsto \sum_{j=0}^\infty
\frac{(f_2^*(t) - f_1^*(t))^j}{j!} \otimes \frac{d^j}{dt^j} \bv;
\]
this is the local analogue of \cite[Proposition~2.2.17]{ber2}.
\end{remark}

\begin{defn} \label{D:h0h1}
For $\calE \in \calM$, write $H^0(\calE)$ for the direct limit of 
$\ker(\nabla)$ on $A[a,1)$ as $a$ approaches 1. Write
$H^1(\calE)$ for the direct limit of $\coker(\nabla)$ on $A[a,1)$ as
$a$ approaches 1. Note that the formation of these commutes with
tensoring over $K$ with a finite Galois extension $K'$, since the trace
map from $K'$ to $K$ commutes with the action of $\nabla$.
Note also that the Yoneda Ext group $\Ext^i(\calE, \calE')$ in 
the category $\calM$ is equal to $H^i(\calE^\dual \otimes \calE')$ for
$i=0, 1$. (For an analogous fact in a more global setting, see
\cite[Proposition~1.1.2]{chiar-lestum}.)
\end{defn}

\begin{example} \label{E:trivial h0h1}
For $\calE = \calO$, we have $H^0(\calO) = K$ generated by 1,
and $H^1(\calO) = K$ generated by $dt/t$.
\end{example}

\subsection{Interlude: semilinear Galois representations}

Before continuing, we gather up a few easy but not necessarily
standard facts about twisted group representations,
which we will use in the construction of the monodromy
representation of a quasi-unipotent $\nabla$-module.

\begin{convention}
All group actions on sets will be right actions.
That is, if the group $G$ acts on the set $S$, we write $s^g$ for
the image of $s \in S$ under $g \in G$, and require the composition
law $s^{gh} = (s^g)^h$.
\end{convention}

\begin{defn}
Let $G$ be a group which acts on a field $F$ (compatibly with
the field operations). A \emph{semilinear 
representation}
of $G$ over $F$ is a finite dimensional $F$-vector space $V$ equipped with
an action of $G$, which is semilinear in the following sense: for $g \in G$,
$c \in F$, and $\bv \in V$, we have $(c\bv)^g = c^g \bv^g$. 
If we define $F\{G\}$ to be the twisted group algebra, in which
$c\{g\} \cdot d\{h\} = c d^{g^{-1}} 
\{gh\}$, then a semilinear representation can
be reinterpreted as a right $F\{G\}$-module.
We say $V$ is \emph{trivial} if it is isomorphic to $F^n$, with the 
$G$-action acting on each copy of $F$ separately, for some nonnegative
integer $n$.
\end{defn}

Maschke's theorem on the complete reducibility of representations
of finite groups goes over to twisted representations as follows.
\begin{lemma}[Maschke property] \label{L:maschke}
Let $G$ be a finite group which acts on a field $F$ of characteristic zero.
Then every semilinear representation of $G$ over $F$ is completely
reducible, i.e., is a direct sum of irreducible twisted representations.
\end{lemma}
\begin{proof}
It suffices to show that if $V$ is indecomposable, then it is also
irreducible. Suppose on the contrary that $V$ is indecomposable, but
$V$ has an irreducible twisted subrepresentation $W$. Choose any projector
$P \in V^\dual \times V$ with image $W$, and put
\[
P' = \frac{1}{\# G} \sum_{g \in G} P^g.
\]
Then $P'$ is again a projector with image $W$, but $P'$ is $G$-invariant.
Hence $1-P'$ is a $G$-invariant projector whose image is a complementary
subrepresentation $W'$ of $W$, contradicting the indecomposability of $V$
and yielding the claim.
\end{proof}

We also have a form of Schur's lemma, with the usual proof.
\begin{lemma}[Schur's lemma] \label{lem:schur}
Let $G$ be a finite group which acts on a field $F$, and let
$V$ and $W$ be irreducible semilinear representations of $G$ over $F$.
Then any $G$-equivariant linear map $f: V \to W$ is either zero or
invertible. In particular, the set of $G$-endomorphisms of
$V$ is a division algebra.
\end{lemma}
\begin{proof}
If $f$ is nonzero, then $\ker(f)$ is a proper
subrepresentation of $V$ and so must vanish, as must $\image(f)$.
\end{proof}

\begin{remark}
Note that semilinear representations may be viewed
as 1-cocycles for $\GL_n(F)$. In particular, if $G = \Gal(F/E)$ for
$F/E$ finite, then $H^1(G, \GL_n(F))$ vanishes,
so any semilinear representation of $G$ over $F$
is trivial.
\end{remark}

\begin{defn}
Let $G \to H$ be a homomorphism of finite groups acting on the field $F$;
this homomorphism induces a homomorphism $F\{G\} \to F\{H\}$ of 
noncommutative $F$-algebras.
Given a semilinear representation $V$ of $G$, we define the
\emph{induced representation} $\Ind^G_H V = V \otimes_{F\{G\}} F\{H\}$;
given a semilinear representation $W$ of $H$, we define the
\emph{restricted representation} $\Res^G_H W = W$ viewed as
a right $F\{G\}$-module via the homomorphism $F\{G\} \to F\{H\}$.
\end{defn}

Since the functors $\Ind^G_H$ and $\Res^G_H$ are left and right
adjoints of each other \cite[Proposition~3.8]{jacobson}, one
obtains the Frobenius reciprocity law as in the linear case.
\begin{lemma}[Frobenius reciprocity] \label{L:Frob rec}
Let $G \to H$ be a homomorphism of finite groups acting on the field $F$,
let $V$ be a semilinear representation of $H$, and let
$W$ be a semilinear representation of $G$. Then there is a natural isomorphism
\[
\Hom_H(V, \Res^G_H W) \cong \Hom_G(\Ind^G_H V, W).
\]
\end{lemma}

\begin{defn}
Let $G$ be a finite group acting on the field $F$.
The \emph{regular representation} of $G$ is the 
semilinear representation corresponding to $F\{G\}$ viewed
as a right module over itself.
\end{defn}

\begin{cor} \label{C:irreducible occurs}
Let $G$ be a finite group acting on the field $F$.
Then any irreducible semilinear representation of $G$ is isomorphic to a
subrepresentation of the regular representation.
\end{cor}
\begin{proof}
Let $V$ be an irreducible semilinear representation of $G$.
By Frobenius reciprocity applied to the trivial group mapping
into $H$, $\Hom_F(F,V) \cong \Hom_G(F\{G\},V)$ is nonzero.
Thus the decomposition of the regular representation into
irreducibles must include a copy of $V$, or else $\Hom_G(F\{G\},V)$
would vanish by Schur's lemma.
\end{proof}

\subsection{Quasi-constant connections}

\begin{defn}
For $\calE \in \calM$, 
we say $\calE$ is \emph{quasi-constant} (or \emph{\'etale})
if there exists $r \in (0,1) \cap \Gamma^*$ such that $\calE$ is defined on
$A[r,1)$, and there exists
a formally \'etale cover $f: X \to A[r,1)$ such that
$f^* \calE$ is spanned by finitely many horizontal sections;
if $R/k((t))$ is the reduction of $X$,
we also say that $\calE$ \emph{is/becomes constant over $X$}
or \emph{over $R$}. Let $\calM^{\qc} = \calM^{\qc}_K$ denote the subcategory
of $\calM$ consisting of quasi-constant objects; it is closed
under formation of direct sums, tensor products, duals, subobjects,
and quotients, though not under formation of extensions.
\end{defn}

One can give a representation-theoretic description of
quasi-constant $\nabla$-modules as follows.
\begin{lemma} \label{L:nabla equivalence}
Let $R/k((t))$ be a finite Galois extension,
and let $f: X \to A[a,1)$ be a formally \'etale cover with
reduction $R$. Let $K'$ be the integral closure of
$K$ in $\Gamma(\calO,X)$. Then the functor $\calE \mapsto H^0(\calE,X)$,
from the category of $\nabla$-modules on $A[a,1)$ which become constant on $R$, 
to the category of semilinear representations of $G = \Gal(R/k((t)))$ in
finite dimensional $K'$-vector spaces, is an equivalence of categories.
\end{lemma}
\begin{proof}
First note that the functor is fully faithful: given
$\calE, \calF$, we have
\[
H^0(\calE, X)^\dual \otimes
H^0(\calF,X) \cong H^0(\calE^\dual \otimes \calF,X) 
\]
since both sides are $K'$-vector spaces of dimension
 $(\rank \calE)(\rank \calF)$
and there is a natural injective map from the left side to the right.
Hence any $G$-equivariant homomorphism between $H^0(\calE,X)$ and
$H^0(\calF,X)$ corresponds to a $G$-equivariant horizontal section of
$\calE^\dual \otimes \calF$ over $X$, and hence to a horizontal
section $\calE^\dual \otimes \calF$ over $A[a,1)$. The
latter corresponds to a morphism from $\calE$ to $\calF$, yielding
the full faithfulness.

Next note that the functor is essentially surjective: put
$\calF = f^* f_* \calO_X$, so that $H^0(\calF,X)$ is the regular representation
of $G$ over $K'$. Let $V$ be an irreducible representation of $G$.
By Corollary~\ref{C:irreducible occurs},
we know that $V$ is isomorphic to
a subrepresentation of $H^0(\calF,X)$. That is, we may choose a $G$-equivariant
projector $P: H^0(\calF,X) \to H^0(\calF,X)$ with image isomorphic to $V$
as a $G$-representation. By the full faithfulness assertion, $P$
comes from a projector $\calF \to \calF$, whose image $\calG$ is a
$\nabla$-module on $A[a,1)$ with $H^0(\calG,X) \cong V$ as
a $G$-representation. Hence the functor is essentially surjective,
and thus an equivalence of categories.
\end{proof}

\subsection{Local duality}

\begin{defn}
For any $\calE \in \calM$, define the
\emph{local duality pairing} on $\calE$ to be the $K$-bilinear pairing
\[
H^0(\calE^\dual) \times H^1(\calE) \to
H^1(\calE^\dual \otimes \calE) \to H^1(\calO);
\]
note again (as in Example~\ref{E:trivial h0h1})
that the latter may be identified with $K$ via the residue map
taking $\sum c_i t^i\,dt$ to $c_{-1}$. We say $\calE$ is 
\emph{dualizable} if the local duality pairing is perfect, i.e.,
if it induces an isomorphism $H^0(\calE^\dual) \cong H^1(\calE)^\dual$.
\end{defn}

\begin{remark} \label{R:duality}
Note that if $\calE = \calE_1 \oplus \calE_2$ in $\calM$, 
then $\calE$ is dualizable
if and only if $\calE_1$ and $\calE_2$ are both dualizable, because the
formation of $H^0$ and $H^1$ commutes with direct sums. Also, if
If $0 \to \calE_1 \to \calE \to \calE_2 \to 0$ is a short exact sequence
in $\calM$, and $\calE_1$ and $\calE_2$ are both dualizable, then
$\calE$ is also dualizable: this follows from applying the snake lemma
to the diagram
\[
\xymatrix{
0 \ar[r] & H^0(\calE_2^\dual) \ar[r] \ar[d] & H^0(\calE^\dual)
\ar[r] \ar[d] & H^0(\calE_1^\dual) \ar[r] \ar[d] & 0 \\
0 \ar[r] & H^1(\calE_2)^\dual \ar[r] & H^1(\calE)^\dual \ar[r] & H^1(\calE_1)^\dual
\ar[r] & 0.
}
\]
\end{remark}

\begin{prop} \label{P:dualizable}
Any $\calE \in \calM^{\qc}$ is dualizable; moreover, any $\calE \in \calM$
which admits a filtration whose successive quotients are quasi-constant
is also dualizable.
\end{prop}
\begin{proof}
The second assertion follows from the first by Remark~\ref{R:duality},
so we will stick to considering $\calE \in \calM^{\qc}$.
Choose a formally \'etale cover $f: X \to A[a,1)$ over which $\calE$ becomes
constant. 
As noted in Definition~\ref{D:h0h1}, the formation of $H^0$ and $H^1$
commutes with tensoring over $K$ with a finite Galois extension of $K$,
so we may enlarge $K$ as needed to ensure that $X \cong A[b,1)$ for some
$b$.

As in the proof of Lemma~\ref{L:nabla equivalence},
we see that $\calE$ is a direct summand of a sum of copies of
$f_* \calO_X$. By
Remark~\ref{R:duality}, to verify that $\calE$ is dualizable, it suffices
to verify that $f_* \calO_X$ is dualizable. But
$H^i(f_* \calO_X) = H^i(\calO_X)$ for $i=0,1$, so it suffices to
verify that $\calO_X$ is dualizable; this in turn follows from the
observation of Example~\ref{E:trivial h0h1} that
$H^0(\calO) = H^1(\calO) = K$.
\end{proof}

\subsection{Unipotent connections}

\begin{defn}
For $\calE \in \calM_K$, we say $\calE$ is \emph{unipotent} if 
$\calE$ admits a filtration whose successive quotients are constant.
Let $\calM^{\unip} = \calM^{\unip}_K$ denote the subcategory of $\calM$
consisting of unipotent objects; it is closed under formation of
direct sums,
tensor products, duals, subobjects, quotients, and extensions.
\end{defn}

Unipotent $\nabla$-modules can be characterized in terms of 
logarithmic connections.
\begin{defn}
Let $\Omega^1_{\log}$ be the coherent sheaf on $A[0,1)$ freely generated by
$\frac{dt}{t}$; we view $\Omega^1$, which is freely generated by $dt$,
as a subsheaf of $\Omega^1_{\log}$. Of course the two coincide away from
$t=0$ (i.e., on $A(0,1)$). A \emph{log-$\nabla$-module} on
$A[0,1)$ is a coherent locally free sheaf $\calE$ equipped with
a connection $\nabla: \calE \to \calE \otimes \Omega^1_{\log}$.
If $\calE$ is a log-$\nabla$-module, then $t \frac{d}{dt}$ acts linearly
on $\Gamma(\calE,A[0,0]) = \calE_0$ (the stalk of $\calE$ at $t=0$), 
which is a finite dimensional vector space
over $K$. We call this vector space equipped with a linear transformation
the \emph{residue} of $\calE$ (or of $\nabla$).
\end{defn}

One then has the following characterization
(compare \cite[Theorem~4.1]{matsuda-katz} and the discussion in
\cite[Chapter~3]{me-part1}).
\begin{prop} \label{P:characterize unip}
The following categories are equivalent.
\begin{enumerate}
\item[(a)]
The category of finite dimensional 
$K$-vector spaces equipped with nilpotent endomorphisms.
\item[(b)]
The category of log-$\nabla$-modules on $A[0,1)$
with nilpotent residue,
whose restrictions to $A(0,1)$ are overconvergent.
\item[(c)]
The category of unipotent $\nabla$-modules on $A(I)$
for any aligned subinterval $I$ of $(0,1)$.
\end{enumerate}
\end{prop}
\begin{proof}
We first exhibit the functor from (a) to (b).
Given a finite dimensional $K$-vector space $V$, put
$M = \calR_{(0,1)} \otimes_K V$ and let $\calE$ be the sheaf
on $A(0,1)$ associated to $M$. We may define a connection on $M$
by declaring that the action of $t \frac{d}{dt}$ on $V$ is via the
given nilpotent endomorphism, then extending via the Leibniz rule.
The resulting log-$\nabla$-module is overconvergent because it is
a successive extension of trivial $\nabla$-modules, and the property
of overconvergence is stable under extensions.

We next exhibit the functor from (b) to (c).
Let $\calE$ be a log-$\nabla$-module on $A[0,1)$ of rank $n$
with nilpotent residue, let $m$ be the index of nilpotency of the
residue map, and let $D$ denote the operator induced by
$t \frac{d}{dt}$. 
Let $P_i$ denote the $i$-th binomial polynomial
\[
P_i(x) = \frac{x(x-1)\cdots(x-i+1)}{i!};
\]
then the $P_i$ form a $\ZZ$-basis for the set of integral-valued
polynomials in $\QQ[x]$.
Let $Q_i$ denote the polynomial
\[
Q_i(x) = x^{m-1} \left( \frac{(1-x)\cdots(i-x)}{i!} \right)^n;
\]
then $Q_{i+1}(x) - Q_i(x)$ is integral-valued of degree $(i+1)n+m-1$
and vanishes at $0,1,\dots, i$, so is an integral linear combination
of $P_{i+1}, \dots, P_{(i+1)n+m-1}$.

By computing on formal power series in $t$
(with which we can formally construct a basis of sections killed by
$D$), we see that for any $b \in [0,1)
\cap \Gamma^*$, $Q_{i+1}(D)-Q_i(D)$
carries any element of $\Gamma(\calE, A[0,b])$
to a multiple of $t^{i+1}$ in the same module. That is,
\[
\frac{1}{t^{i+1}}(Q_{i+1} - Q_i)(D) 
\]
is a well-defined operator on $\calE$.
As we saw above, $\frac{1}{t^{i+1}} (Q_{i+1}-Q_i)(D)$
is an integer 
linear combination of
\[
\frac{1}{t^{i+1}} P_l(D) \qquad (l = i+1, \dots, (i+1)n+m-1),
\]
and hence is a $\calR_{[0,1)}^{\inte}$-linear 
combination of the $\frac{1}{t^l} P_l(x)$
for $l=i+1, \dots, (i+1)n+m-1$.

However,
\[
\frac{1}{t^l}
P_l(D) = 
\frac{1}{l!} \frac{d^l}{d t^l};
\]
by the overconvergence condition,
for any $\eta \in (0,1)$, there exists $b \in [0,1) \cap \Gamma^*$
such that for any $c \in (b,1) \cap \Gamma^*$ and
any $\bv \in \Gamma(\calE, A[b,c])$,
the sequence 
\[
\left\{\frac{1}{i!} \frac{d^i}{dt^i} \bv \right\}_{i=1}^\infty
\]
is $\eta$-null.
It follows that the sequence
\begin{equation} \label{eq:unipotent sequence}
\{t^{-i-1} (Q_{i+1} - Q_i)(D)\bv\}
\end{equation}
 is also $\eta$-null over $A[b,c]$.

Now choose $\bv \in \Gamma(\calE, A[0,c])$
with $c \in (\eta,1) \cap \Gamma^*$.
By Lemma~\ref{L:Hadamard}, the sequence \eqref{eq:unipotent sequence}
is $\eta$-null not just over
$A[b,c]$, but also over $A[0,c]$.
In particular, the sequence \eqref{eq:unipotent sequence}
is $\eta$-null over $A[0,\eta]$, and so the sequence
$\{(Q_{i+1} - Q_i)(D) \bv\}$ is 
$1$-null over $A[0,\eta]$.
That is, the limit
\[
f(\bv) = \lim_{i \to \infty} Q_i(D) \bv
\]
exists in $\Gamma(\calE, A[0,\eta])$.

Again from the formal power series computation, we see that
$f$ acts as the $(m-1)$-st power of the residue map modulo $t$, and
that $Df(\bv) = 0$ for all $\bv$. We can find some $\bv \in 
\Gamma(\calE, A[0,c])$ whose fibre at zero is not killed by the
$(m-1)$-st power of the residue map; then $f(\bv) \neq 0$ but
$Df(\bv) = 0$. That is, the log-$\nabla$-submodule of $\calE$ spanned by
the kernel of $\nabla$ on $A[0,c]$ is nonzero. The rank of that submodule
does not increase as $c$ increases, so it must be constant for $c \in (0,1)
\cap \Gamma^*$ sufficiently closed to 1. Thus these submodules fit together
to yield a constant $\nabla$-submodule of $\calE$ on $A[0,1)$; quotienting
by this submodule and repeating the argument, we obtain the desired
unipotent filtration.

Finally, we note that the functor from (c) to (a) is straightforward:
it suffices to verify that on each open subinterval of $I$,
$\calE$ is spanned by sections killed by $D^n$ with $n = \rank(\calE)$. 
That in turn
follows by induction on rank, using the fact that for $I$ open,
the cokernel of $t \frac{d}{dt}$ on $\calR_I$ is generated over $K$ by 1.
\end{proof}

\begin{remark}
This discussion goes over to higher-dimensional polydiscs;
see \cite[Chapter~3]{me-part1}.
We recall also a remark from \cite[Chapter~3]{me-part1}:
the application of 
Lemma~\ref{L:Hadamard} must be to a sequence without poles,
which necessitates the introduction of the sequence
\eqref{eq:unipotent sequence} in lieu of working directly
on the sequence $\{\frac{1}{i!} \frac{d^i}{dt^i} \bv\}$.
\end{remark}

Proposition~\ref{P:characterize unip} reduces most questions about
unipotent $\nabla$-modules to linear algebra, as in the following
special case of \cite[Lemma~7.6]{matsuda-katz}.
\begin{lemma} \label{L:unip h0h1}
If $\calU \in \calM^{\unip}$ is nonzero and indecomposable, then
$\dim_K \Ext^i(\calU) = 1$ for $i=0,1$.
\end{lemma}
\begin{proof}
Note that any element of $\calM^{\unip}$ is dualizable by
Proposition~\ref{P:dualizable}, so it suffices to prove the
claim for $i=0$. By Proposition~\ref{P:characterize unip},
we may translate the claim into linear algebraic terms: the result
is simply the fact that a nilpotent linear transformation on a finite
dimensional vector space is indecomposable if and only if it can
be written as a single Jordan block.
\end{proof}

We also need the following result on the interaction between quasi-constant
and unipotent $\nabla$-modules.
\begin{lemma} \label{L:h0 qc}
For $\calP \in \calM^{\qc}$ and $\calU \in \calM^{\unip}$, the natural map
$H^0(\calP) \otimes_K H^0(\calU) \to H^0(\calP \otimes \calU)$ is an isomorphism
of $K$-vector spaces.
\end{lemma}
\begin{proof}
If $\calU = \calO$, this is obvious; otherwise, we proceed by induction
on $\rank(\calU)$. Choose a short exact sequence
$0 \to \calO \to \calU \to \calU_1 \to 0$, and consider the diagram
\begin{equation} \label{eq:snake diag0}
\xymatrix{
0 \ar[r] & H^0(\calP) \ar[d] \ar[r] & H^0(\calP) \otimes H^0(\calU)
\ar[d] \ar[r] & H^0(\calP) \otimes H^0(\calU_1) \ar[r] \ar[d] & H^0(\calP) 
\otimes H^1(\calO)
\ar[d] \\
0 \ar[r] & H^0(\calP) \ar[r] & H^0(\calP \otimes \calU) \ar[r] &
H^0(\calP \otimes \calU_1) \ar[r] & H^1(\calP)
}
\end{equation}
in which the first
row is obtained by applying the snake lemma to the diagram
\begin{equation} \label{eq:snake diag}
\xymatrix{
0 \ar[r] & \calO \ar^d[d] \ar[r] & \calU \ar^\nabla[d] \ar[r] &
\calU_1 \ar^\nabla[d] \ar[r] & 0 \\
0 \ar[r] & \Omega^1 \ar[r] & \calU \otimes \Omega^1
\ar[r] & \calU_1 \otimes \Omega^1 \ar[r] & 0
}
\end{equation}
and then tensoring with $H^0(\calP)$, whereas the second row is obtained
by first tensoring \eqref{eq:snake diag} with $\calP$ and then applying the
snake lemma. In \eqref{eq:snake diag0}, the first vertical arrow is 
visibly an isomorphism, and the third is an isomorphism by the induction
hypothesis. As for the fourth arrow, note that $H^1(\calP) \cong
H^0(\calP^\dual)^\dual$ by Proposition~\ref{P:dualizable},
whereas by Lemma~\ref{L:nabla equivalence}, forming $H^0$ of a
quasi-constant $\nabla$-module commutes with taking duals
(since a semilinear representation is trivial if and only if its dual
is trivial). Hence the fourth vertical arrow is also an isomorphism;
by the five lemma, the arrow $H^0(\calP) \otimes H^0(\calU)
\to H^0(\calP \otimes \calU)$ is an isomorphism, as desired.
\end{proof}

\subsection{Quasi-unipotent connections}

\begin{defn}
For $\calE \in \calM$,
we say $\calE$ is \emph{quasi-unipotent} if $\calE$ admits a
filtration whose successive quotients are
quasi-constant; if each successive quotient becomes constant
over $X$ or over $R$, we say
we say $\calE$ \emph{is/becomes unipotent over $X$} or
\emph{over $R$}.
Let $\calM^{\qu} = \calM^{\qu}_K$ denote the subcategory of $\calM_K$
consisting of quasi-unipotent objects; it is closed under formation of
direct sums,
tensor products, duals, subobjects, quotients, and extensions.
Note that any element of $\calM^{\qu}$ is dualizable by
Proposition~\ref{P:dualizable}.
\end{defn}
\begin{remark}
The definition of quasi-unipotence in \cite{matsuda-katz} is slightly different:
there it is required that $\calE$ admit a unipotent filtration over
$X$ for some formally \'etale cover $f: X \to A[r,1)$. 
This clearly follows from the existence of a filtration over $A[r,1)$
whose successive quotients are quasi-constant, but in fact the reverse is
also true: if $\calE$ admits a unipotent filtration over $X$, then there
is a unique such filtration of shortest length, namely the one whose first
step is spanned by the full kernel of $\nabla$ over $X$, and so on.
By Galois descent, this filtration descends to $A[r,1)$.
\end{remark}

Quasi-unipotent $\nabla$-modules can be canonically decomposed into
``isotypical'' pieces.
\begin{defn}
Suppose $\calF \in \calM$ is irreducible. We say that $\calE \in \calM$
is \emph{$\calF$-typical} if $\calE$ admits a filtration whose successive
quotients are all isomorphic to $\calF$.
\end{defn}
In this language, we have the following partial analogue of
\cite[Lemma~7.6]{matsuda-katz}.
\begin{lemma} \label{L:no hom}
Suppose that $\calE_j$ is $\calF_j$-typical for $j=1,2$, and that
$\calF_1 \not\cong \calF_2$. Then $\Ext^i(\calE_1, \calE_2) = 0$
for $i=0,1$ (in the category $\calM$).
\end{lemma}
\begin{proof}
As noted earlier (see Definition~\ref{D:h0h1}), 
$\Ext^i(\calE_1, \calE_2) \cong H^i(\calE_1^\dual \otimes \calE_2)$.
By Proposition~\ref{P:dualizable}, it suffices to check the claim for
$i=0$, in which case it follows by induction on rank plus Schur's lemma.
\end{proof}

\begin{prop} \label{P:isotypical decomp}
Each $\calE \in \calM$ admits a unique (up to reordering) decomposition
$\calE_1 \oplus \cdots \oplus \calE_m$, in which each
$\calE_i$ is $\calF_i$-typical for some irreducible $\calF_i \in \calM$,
and no two of the $\calF_i$ are isomorphic.
\end{prop}
\begin{proof}
For existence, we proceed by induction on $\rank(\calE)$. If
$\calE$ is irreducible, there is nothing to check; otherwise,
choose an exact sequence $0 \to \calF \to \calE \to \calE' \to 0$
with $\calF$ irreducible. By the induction hypothesis,
$\calE'$ admits a decomposition 
$\calE'_1 \oplus \cdots \oplus \calE'_m$ of the desired form.
Let $\calE_i$ be the inverse image of $\calE'_i$ in $\calE$;
for all but possibly one index $i$, $\calE'_i$ is $\calF_i$-isotypical
for some $\calF_i \not\cong \calF$, and so the exact sequence
$0 \to \calF \to \calE_i \to \calE'_i \to 0$ splits by
Lemma~\ref{L:no hom}. We may thus decompose $\calE$ in the desired form.

For uniqueness, note that if
$\calE_1 \oplus \cdots \oplus \calE_m \cong \calE'_1 \oplus \cdots \oplus
\calE'_m$, then the isomorphism must carry $\calE_1$ into an
$\calF_1$-typical summand by Lemma~\ref{L:no hom}, and so on.
\end{proof}

One can further refine Proposition~\ref{P:isotypical decomp}; the result is
essentially \cite[Theorem~7.8]{matsuda-katz}.
\begin{theorem}[Matsuda] \label{T:Matsuda}
Every $\calE \in \calM_K^{\qu}$ admits
a canonical decomposition as a direct sum
$\sum_i \calF_i \otimes \calU_i$, where each
$\calF_i$ is quasi-constant and irreducible, each $\calU_i$ is unipotent,
and no two of the $\calF_i$ are isomorphic.
\end{theorem}
\begin{proof}
By virtue of Proposition~\ref{P:isotypical decomp}, it suffices to
check that for $\calF \in \calM^{\qc}$ irreducible, any $\calF$-typical
element $\calE$ of $\calM$ has the form $\calF \otimes \calU$ for some
$\calU \in \calM^{\unip}$. 
We check this by induction on $\rank(\calE)$; we may of course assume
$\calE$ is indecomposable.

If $\calE$ is irreducible, the claim is clear; otherwise,
construct a short
exact sequence $0 \to \calE_1 \to \calE \to \calF \to 0$, and note that
$\calE_1$ is $\calF$-typical.
By the induction hypothesis,
we have $\calE_1 \cong \calF \otimes \calU_1$ for some
$\calU_1 \in \calM^{\unip}$, necessarily indecomposable.
By Proposition~\ref{P:dualizable} and Lemma~\ref{L:h0 qc},
\begin{align*}
\Ext^1(\calF, \calE_1) &= 
H^1(\calF^\dual \otimes \calE_1) \\
&= H^0(\calE_1^\dual \otimes \calF)^\dual \\
&= H^0(\calF^\dual \otimes \calF \otimes \calU_1^\dual)^\dual \\
&= H^0(\calF^\dual \otimes \calF)^\dual \otimes H^0(\calU_1^\dual)^\dual \\
&= H^0(\calF^\dual \otimes \calF)^\dual \otimes H^1(\calU_1)  \\
&= H^0(\calF^\dual \otimes \calF)^\dual \otimes \Ext^1(\calO, \calU_1).
\end{align*}
Let $\calU_1 = \oplus_i \calU_{1,i}$ be the decomposition of
$\calU_1$ into indecomposables.
By Lemma~\ref{L:unip h0h1}, for each $i$, $\dim_K
\Ext^1(\calO,\calU_{1,i}) = 1$.
We may thus adjust each component of $\calE_1 \cong \oplus_i (\calF \otimes
\calU_{1,i})$ by an automorphism of $\calF$ to ensure that
the element of $\Ext^1(\calF,\calE_1)$ corresponding to the exact sequence
$0 \to \calE_1 \to \calE \to \calF \to 0$ maps to an element of
$H^0(\calF^\dual \otimes \calF)^\dual \otimes \Ext^1(\calO, \calU_1)$
in the image of the map $1 \otimes \id$.
In this case, the exact sequence 
$0 \to \calE_1 \to \calE \to \calF \to 0$ is obtained
from some exact sequence $0 \to \calU_1 \to \calU \to \calO
\to 0$ by tensoring with $\calF$; in particular,
$\calE \cong \calU \otimes \calF$, as desired.
\end{proof}

\begin{remark}
One can deduce from Theorem~\ref{T:Matsuda} that
every quasi-unipotent $\nabla$-module is overconvergent,
by checking in the unipotent and quasi-constant cases.
In the latter case, one can construct a so-called ``unit-root'' Frobenius
structure as in \cite[Lemma~5.3]{matsuda-katz} and then apply
Lemma~\ref{L:Frobenius convergent} below; alternatively, one can
check that the overconvergence property descends down formally
\'etale covers.
\end{remark}

\subsection{Monodromy representations}

Theorem~\ref{T:Matsuda} allows us to describe the
category of quasi-unipotent $\nabla$-modules in
represen\-tation-theoretic terms, as follows.
\begin{defn}
Put $G^{\log} = 
G_K^{\log} = \Gal(k((t))^{\sep}/k((t))) \times K$, where the second
factor represents the additive group. A semilinear representation
of $G^{\log}$ on a finite dimensional $K^{\unr}$-vector space is
\emph{permissible} if its restriction to some open subgroup
of $\Gal(k((t))^{\sep}/k((t)))$ is trivial (in the sense of
being isomorphic to a product of copies of $K^{\unr}$), and
its restriction to $K$ is algebraic (and hence necessarily unipotent).
\end{defn}

\begin{defn}
Define the sheaf $\calO_{\log}$ on $A[a,1)$ by the formula
\[
\Gamma(\calO_{\log}, A[b,c]) = \Gamma(\calO, A[b,c])[\log(t)],
\]
where $\log(t)$ is an indeterminate; define the sheaf
$\calO_{\log}$ on a finite \'etale cover of $A[a,1)$
as the pullback from $A[a,1)$.
Extend $d$ to a $K$-derivation $\calO_{\log} \to \Omega^1_{\log}$
by setting $d(\log(t)) = \frac{dt}{t}$.
For $\calE \in \calM$ and $f: X \to A[a,1)$ a formally
\'etale cover, let $H^0_{\log}(\calE,X)$ denote the kernel
of $\nabla$ on $\Gamma(\calO_{\log}, X)$.
\end{defn}
\begin{remark}
Beware that one has a canonical isomorphism $\tau^* \calO_{\log} \to 
\calO_{\log}$ (i.e., an action of $\tau$ on $\calO^{\log}$) not for an
arbitrary automorphism
$\tau$ of $A[a,1)$, but only for those $\tau$ satisfying
$|\tau^*(t) - t|_r < 1$ for $r \in [a,1)$; the point is that
these are the $\tau$ for which the series defining
$\log (\tau^*(t)/t)$ converges in $\Gamma(\gotho, A[a,1))$.
\end{remark}

\begin{lemma}
Let $f: X \to A[a,1)$ be a formally \'etale cover over 
which $\calE \in \calM^{\qu}_K$ becomes unipotent, and let
$K'$ be the integral closure of $K$ in $\Gamma(\calO,X)$. Then
$\dim_{K'} H^0_{\log}(\calE, X) = \rank(\calE)$. 
\end{lemma}
\begin{proof}
By Theorem~\ref{T:Matsuda}, it is enough to verify this
for $\calE$ quasi-constant, in which case it is evident,
and for $\calE$ unipotent, in which case it follows from
Proposition~\ref{P:characterize unip}.
\end{proof}

\begin{defn}
For $\calE \in \calM^{\qu}$, 
let $f: X \to A[a,1)$ be a formally \'etale cover over
which $\calE$ becomes unipotent, and let $K'$ be the integral
closure of $K$ in $\Gamma(\calO,X)$. Then we may view
\[
H^0_{\log}(\calE, X) \otimes_{K'} K^{\unr}
\]
as a semilinear representation of
$G^{\log}$, where the action of $c \in K$ is by the map
$\log(x) \mapsto \log(x) + c$. This representation does not
depend on the choice of $X$; we call it the 
\emph{monodromy representation} of $\calE$.
\end{defn}

We can now give our representation-theoretic description
of $\calM^{\qu}$.
\begin{theorem} \label{T:monodromy rep}
The monodromy representation, viewed as a functor from
$\calM^{\qu}$ to the category of permissible semilinear representations
of $G^{\log}$, is an equivalence of categories.
\end{theorem}
\begin{proof}
The full faithfulness of the functor follows as in the proof of
Lemma~\ref{L:nabla equivalence}. As for the essential surjectivity,
note that any indecomposable representation of the product group
$G^{\log}$ is the tensor product of a permissible representation of each
factor. Each permissible representation of 
$\Gal(k((t))^{\sep}/k((t)))$ occurs in $\calM^{\qu}$
by Lemma~\ref{L:nabla equivalence}; each permissible representation
of $K$ corresponds to a nilpotent linear transformation, and hence
to a unipotent $\nabla$-module by 
Proposition~\ref{P:characterize unip}. Thus the functor is
essentially surjective, completing the proof.
\end{proof}
\begin{remark}
For $k$ algebraically closed (and $K$ discretely valued, though this is
less critical), Theorem~\ref{T:monodromy rep} is essentially due to
Andr\'e \cite[Th\'eor\`eme~7.1.1]{andre}.
\end{remark}

\section{Radius of convergence and ramification}
\label{sec:cm}

In this chapter, we recall a method of Christol and Mebkhout
that attaches numerical invariants to a quasi-unipotent $\nabla$-module,
and a theorem of Matsuda that relates these to the ramification breaks.
Note that our derivation of Matsuda's theorem is direct, and does not
depend \emph{per se} on the theory of Christol-Mebkhout.

We retain all notations and conventions from the preceding chapters,
notably Convention~\ref{conv:mono}.

\subsection{Generic radius of convergence}

We first recall the notion of ``generic radius of convergence'' from
Christol-Dwork \cite{christol-dwork}.
\begin{defn} \label{D:operator norm}
The \emph{operator norm} $|T|_{\op}$
of a continuous operator $T$ on a normed space
$V$ is equal to the smallest nonnegative real number $c$ such
that $|T\bv| \leq c|\bv|$ for all $\bv \in V$. 
For instance,
if $V = \calR_I$ equipped with the norm $|\cdot|_r$ for
some $r \in I$, and $T = \frac{d}{dt}$, then $|T\bv|_r \leq r^{-1} |\bv|_r$
with equality for $\bv = t$, so the operator norm equals $r^{-1}$.
\end{defn}

\begin{defn}
The \emph{spectral norm} $|T|_{\spec}$
of a continuous operator $T$ on a normed space
$V$ is defined as
\[
|T|_{\spec} = \limsup_{s \to \infty} |T^s|_{\op}^{1/s}.
\]
For instance,
if $V = \calR_I$ equipped with the norm $|\cdot|_r$ for
some $r \in I$, and $T = \frac{d}{dt}$, then $|T^s \bv|_r \leq r^{-s} |s!| |\bv|_r$
with equality for $\bv = t^s$, so 
\[
\left|\frac{d}{dt} \right|_{\spec} = \limsup_{s \to \infty} 
\,(r^{-s} |s!|)^{1/s} = p^{-1/(p-1)} r^{-1}.
\]
\end{defn}

\begin{defn}
Let $I$ be a closed aligned subinterval of $[0,1)$,
and let $\calE$ be a $\nabla$-module over $A(I)$;
put $M = \Gamma(\calE,A_K(I))$. Since $\calR_I$ is a
principal ideal domain (Proposition~\ref{P:poly generator}),
$M$ is free over $\calR_I$; let $\be_1, \dots, \be_n$ be a basis.
Let $T$ denote the operator induced by $\frac{d}{dt}$ on $M$,
and define the matrix $D_s$ by
\[
T^s \be_j = \sum_i (D_s)_{ij} \be_i.
\]
For $\rho \in I$, define the \emph{$\rho$-spectral norm} of $T$ 
(with respect to $\be_1, \dots, \be_n$) to be
\[
|T|_{\rho,\spec} = \min \{p^{-1/(p-1)} \rho^{-1}, \limsup_{s \to \infty}
\max_{i,j} |(D_s)_{ij}|^{1/s} \};
\]
then a short calculation \cite[Proposition~1.2]{christol-dwork}
shows that the $\rho$-spectral norm is independent of the choice of
basis.
\end{defn}

\begin{defn}
Let $I$ be an aligned subinterval of $[0,1)$, and let
$\calE$ be a $\nabla$-module on $A(I)$.
Define the function $\rho \mapsto R(\calE,\rho)$ on $I$ as follows:
\[
R(\calE,\rho) = p^{-1/(p-1)} |T|^{-1}_{\rho,\spec},
\]
where $T$ is the operator induced by $\frac{d}{dt}$ on
$\Gamma(\calE,J)$ for any closed aligned subinterval $J$ of $I$
containing $\rho$. (Note that the choice of $J$ does not matter: one
can compare the spectral norms obtained for one $J$ and a strictly
smaller $J$ by using the same basis to compute them. Indeed,
one could even take $J = [\rho,\rho]$.)
\end{defn}

\begin{remark} \label{R:continuous}
The function $\rho \mapsto R(\calE,\rho)$ is log-concave
\cite[Proposition~2.3]{christol-dwork}, so in particular it is
continuous on the interior of $I$. It also turns out to be
continuous at the endpoints of $I$, but this is subtler
\cite[Th\'eor\`eme~2.3]{christol-dwork}.
One can make similar continuity statements using the Berkovich
foundations of rigid analytic geometry, but they carry somewhat more
content; see \cite{balda-divizio}.
\end{remark}

\begin{remark}
In terms of the function $R$, the condition that $\calE$ be
overconvergent in our sense
is precisely that $\lim_{\rho \to 1} R(\calE,\rho)\rho^{-1} = 1$;
this is literally the condition that $\calE$ be ``soluble at 1''
in the terminology of \cite[4.1-1]{cm3}.
\end{remark}

\subsection{Generic points}

The interpretation of $R(\calE,\rho)$ as a ``generic radius
of convergence'' is as follows.
\begin{defn}
Let $\CC$ be an algebraically closed extension of $K$ complete
for an absolute value extending $|\cdot|$ on $K$, whose
residue field is transcendental over that of $K$. 
For $\rho \in (0,1) \cap \Gamma^*$, a \emph{generic point of radius $\rho$} is an element
$t_{\rho} \in \CC$ such that $|t_{\rho} - x| = \rho$ for any $x \in
K^{\alg}$ with $|x| \leq \rho$. 
(Since this includes $x=0$, we have in particular that $|t_\rho| = \rho$.) 
\end{defn}

\begin{lemma} \label{L:generic point}
Let $I$ be an aligned subinterval of $(0,1)$, and 
suppose $\rho \in I \cap \Gamma^*$.
For any $x = \sum c_i t^i \in \calR_I$ and
any generic point $t_\rho \in \CC$ of radius $\rho$,
\[
|x|_\rho = \left|\sum c_i t_\rho^i\right|.
\]
\end{lemma}
\begin{proof}
Since the supremum $\sup \{|c_i| \rho^i\}$
defining $|x|_\rho$
is achieved by at least one $i$, it suffices
to check the equality for $x \in K[t,t^{-1}]$. Moreover,
multiplying $x$ by $t$ multiplies
both sides of the proposed equality by $\rho$, so we may reduce
to the case where $x \in K[t]$.

Factor $x = c_n \prod_{j=1}^n (t - z_j)$ with each $z_j \in K^{\alg}$.
Then $|t_\rho - z_j| = \max\{\rho, |z_j|\} = |t-z_j|_\rho$
for each $j$, because $t_\rho$ is a generic point. Hence
\begin{align*}
|x|_\rho &\geq \left|\sum c_i t_\rho^i\right| \\
&= |c_n| \prod_{j=1}^n |t_\rho - z_j| \\
&= |c_n| \prod_{j=1}^n |t - z_j|_\rho \\
&\geq |x|_\rho,
\end{align*}
yielding the desired equality.
\end{proof}

\begin{prop} \label{P:generic points}
Let $I$ be an aligned subinterval of $[0,1)$, and let
$\calE$ be a $\nabla$-module on $A(I)$.
For $\rho \in I \cap \Gamma^*$, $t_\rho \in \CC$ a generic
point of radius $\rho$, and $r \in \RR$, the following are equivalent.
\begin{enumerate}
\item[(a)]
$|R(\calE, \rho)| \geq r$;
\item[(b)]
for any $z \in \CC$ with $|z| = \rho$, 
$\calE$ admits a basis of horizontal sections in the disc
$|t-z| < r$;
\item[(c)]
$\calE$ admits a basis of horizontal sections in the disc
$|t-t_\rho| < r$.
\end{enumerate}
\end{prop}
\begin{proof}
Note that for $\bv \in \calE$, the sum
\[
f(\bv) = \sum_{s=0}^\infty \frac{(t-z)^s}{s!} \frac{d^s}{dt^s} \bv,
\]
if convergent, yields a horizontal section of $\calE$ on a disc
centered at $z$.
Thus given (a), we may evaluate $f(\bv)$ on any basis of $\calE$ to
obtain horizontal sections in the disc $|t-z|<r$; that is,
(a) implies (b). Also, (b) implies (c) trivially. On the other hand,
given (c), one deduces $|R(\calE, \rho)| \geq r$ from
Lemma~\ref{L:generic point}.
\end{proof}
\begin{cor} \label{C:generic radius auto}
Given an automorphism $\phi$ of $A[0,1)$, there exists
$\rho_0 \in [0,1) \cap \Gamma^*$ such that $R(\calE,\rho) = 
R(\phi^* \calE, \rho)$ for $\rho \in [\rho,1) \cap \Gamma^*$.
\end{cor}
\begin{proof}
This follows from Proposition~\ref{P:generic points} and the fact
that $\phi$ preserves $A[\rho_0,1)$ for $\rho_0$ sufficiently
close to 1.
\end{proof}

\subsection{Formally \'etale covers and generic points}

We will need to make a few calculations concerning the way discs around generic
points transform under formally \'etale covers.

\begin{lemma} \label{L:generic tame}
Let $n$ be an integer relatively prime to $p$.
For $\rho \in (0,\infty) \cap \Gamma^*$, choose $t_\rho \in \CC$ with
$|t_\rho| = \rho$, and choose an $n$-th root
$t_\rho^{1/n}$ of $t_\rho$.
Then for any $r \in (0, 1] \cap \Gamma^*$,
the map $z \mapsto z^n$ induces an isomorphism between the discs 
\[
|t - t_\rho^{1/n}| < r \rho^{1/n}
\qquad \mbox{and} \qquad
|t - t_\rho| < r \rho.
\]
\end{lemma}
\begin{proof}
If $|z - t_\rho^{1/n}| < \rho^{1/n}$, then
\begin{align*}
|z^n - t_\rho| 
&= \rho |(1 + (z-t^{1/n}_\rho)/t_\rho^{1/n})^n - 1| \\
&= \rho \left|\sum_{i=1}^\infty \binom{n}{i} ((z-t_\rho^{1/n})/t_\rho^{1/n})^i\right| \\
&= \rho |(z-t_\rho^{1/n})/t_\rho^{1/n}|
\end{align*}
since $n$ is relatively prime to $p$.
Thus the map $z \mapsto z^n$ induces a map from
the disc $|t - t_\rho^{1/n}| < r \rho^{1/n}$ into the disc
$|t - t_\rho| < r \rho$.

We define the inverse map by the binomial series
\[
z \mapsto t_\rho^{1/n} \sum_{i=0}^\infty 
\binom{1/n}{i} \left ( \frac{z}{t_\rho} - 1 \right)^i.
\]
If $|z - t_\rho| < r \rho$, then $|(z/t_\rho) - 1| < r$, so the series
converges to a value in the disc $|t - t_\rho^{1/n}| < r \rho^{1/n}$.
This yields the desired result.
\end{proof}

The situation is a bit different when $n$ is not coprime to $p$;
it will be enough for us to consider $n=p$.
\begin{lemma} \label{L:generic Frob}
For $\rho \in (0,\infty) \cap \Gamma^*$, choose $t_\rho \in \CC$ with
$|t_\rho| = \rho$, and choose a $p$-th root
$t_\rho^{1/p}$ of $t_\rho$.
Then for any $r \in (0, p^{-p/(p-1)}] \cap \Gamma^*$,
the map $z \mapsto z^p$ induces an isomorphism between the discs 
\[
|t - t_\rho^{1/p}| < r p \rho^{1/p}
\qquad \mbox{and} \qquad
|t - t_\rho| < r \rho;
\]
for $r \in (p^{-p/(p-1)}, 1) \cap \Gamma^*$,
the disc
\[
|t - t_\rho^{1/p}| < r \rho^{1/p}
\qquad \mbox{is carried into the disc} \qquad
|t-t_\rho| < r^p \rho.
\]
\end{lemma}
\begin{proof}
If $|z - t_\rho^{1/p}| < p^{-1/(p-1)} \rho^{1/p}$, then
\begin{align*}
|z^p - t_\rho| 
&= \rho |(1 + (z-t^{1/p}_\rho)/t_\rho^{1/p})^p - 1| \\
&= \rho |\sum_{i=1}^p \binom{p}{i} ((z-t_\rho^{1/p})/t_\rho^{1/p})^i| \\
&= \rho |p (z-t_\rho^{1/p})/t_\rho^{1/p}|.
\end{align*}
Thus for $r \leq p^{-p/(p-1)}$,
the map $x \mapsto x^p$ induces a map from
the disc $|t - t_\rho^{1/p}| < r p \rho^{1/p}$ into the disc
$|t - t_\rho| < r \rho$.
For $r > p^{-p/(p-1)}$,
the argument breaks down because the dominant term in the binomial
expansion is no longer the first term, but the last term;
however, we still conclude that the disc
$|t - t_\rho^{1/p}| < r \rho^{1/p}$
maps into the disc $|t-t_\rho| < r^p \rho$.

For $r \leq p^{-p/(p-1)}$,
we again define the inverse map by the binomial series
\[
z \mapsto t_\rho^{1/p} \sum_{i=0}^\infty \binom{1/p}{i} 
\left ( \frac{z}{t_\rho} - 1 \right)^i.
\]
This time, however, if 
$|z - t_\rho| < r \rho$, then $|(z/t_\rho) - 1| < r < p^{-p/(p-1)}$, 
so the series
converges to a value in the disc $|t - t_\rho^{1/p}| < r p \rho^{1/p}$.
This yields the desired result.
\end{proof}

Finally, we consider a cover corresponding to a wildly ramified cover
of $\Spec k((t))$.
\begin{lemma} \label{L:generic wild}
Suppose that $K$ contains an element $\pi$ with $\pi^{p-1} = -p$.
For $\rho \in (p^{-p/(p-1)},1) \cap \Gamma^*$, choose $t_\rho \in \CC$ with
$|t_\rho| = \rho$, and
choose $u_\rho \in \CC$ such that 
$(1 + \pi u_\rho)^p = 1 + p \pi t_\rho^{-1}$.
Then $|u_\rho| = \rho^{-1/p}$,
and for $r \in [0, \rho]$, the ring inclusion
\[
\calR_{[\rho,1)} \to 
\calR_{[\rho,1)}[u]/((1 + \pi u)^p - (1 + p \pi t^{-1}))
\]
induces an isomorphism between the discs
\[
|u^{-1} - u_\rho^{-1}| < r \rho^{(2-p)/p}
\qquad \mbox{and} \qquad
|t - t_\rho| < r \rho.
\]
\end{lemma}
\begin{proof}
The polynomial defining $u_\rho$ may be rewritten
\[
0 = 1 + \sum_{i=1}^p \binom{p}{i} (u_\rho t_\rho^{1/p})^i (t_\rho^{1/p} \pi^{-1})^{p-i},
\]
and the resulting polynomial in $u_\rho t_\rho^{1/p}$ has vertices corresponding
to the terms of degree 0 and $p$. Hence $|u_\rho t_\rho^{1/p}| = 1$,
i.e., $|u_\rho| = \rho^{-1/p}$.

By Lemma~\ref{L:generic tame}, the discs
$|t - t_\rho| < r \rho$ and $|t^{-1} - t_\rho^{-1}| < r \rho^{-1}$
are identical, and the discs
$|u^{-1} - u_\rho^{-1}| < r \rho^{(1 -(p-1))/p}$
and $|u - u_\rho| < r \rho^{(-1 -(p-1))/p}$ are identical.
By Lemma~\ref{L:generic Frob} (and the fact that
$|1+\pi u| = 1$ because $\rho > p^{-p/(p-1)}$), the ring inclusion induces
an isomorphism between the discs
\[
|(1+p \pi t^{-1}) - (1+p \pi t_\rho^{-1})| < r p^{-p/(p-1)} \rho^{-1}
\]
and
\[
|(1 + \pi u) - (1 + \pi u_\rho)| < r p^{-1/(p-1)} \rho^{-1}.
\]
This yields the desired result.
\end{proof}

\subsection{Highest breaks and radii of convergence}

We now recall a relationship between generic radii of 
convergence of a quasi-unipotent $\nabla$-module
and the breaks of its monodromy representation
(Theorem~\ref{T:breaks}).
For more on the history of this relationship, see
Remark~\ref{R:breaks}.

\begin{defn} \label{D:highest break}
Given $\calE \in \calM$ and $\beta \in \RR_{\geq 0}$, we say that
$\calE$ has \emph{highest break $\leq \beta$} (resp.\ $\geq \beta$) if there exists
$\rho_0 \in (0,1)$ such that 
\[
R(\calE, \rho) \geq \rho^{\beta + 1}
\quad \mbox{(resp.\ $R(\calE,\rho) \leq \rho^{\beta+1}$)} 
 \qquad \mbox{for $\rho \in (\rho_0,1) \cap
\Gamma^*$}.
\]
We say $\calE$ has highest break $\beta$ if it has highest break
$\leq \beta$ and $\geq \beta$.
Note that this definition is stable under pullback by automorphisms
of the disc, thanks to Corollary~\ref{C:generic radius auto}.
\end{defn}
\begin{remark}
Note that Definition~\ref{D:highest break} does not guarantee that
$\calE$ has a highest break. However, if $\calE$ is overconvergent and $K$
is spherically complete, then
by \cite[Th\'eor\`eme 4.2-1]{cm3}
and \cite[Th\'eor\`eme 2.1-2]{cm4}
(and an application of Theorem~\ref{T:coherent} to convert $\calE$
into a module over the Robba ring), 
$\calE$ has highest break $\beta$ for some $\beta \in \QQ \cap [0, \infty)$.
We will not explicitly use this result; however, note that it is 
implicit in the proofs of the $p$-adic local monodromy theorem
given in \cite{andre} and \cite{mebkhout}. (On the other hand,
it plays no role in the proof given in \cite{me-local}.)
\end{remark}

\begin{remark}
Christol and Mebkhout use the term ``plus grande pente'', which
translates as ``greatest slope'', for what we are calling the
``highest break''. We have avoided the term ``slope'', which
Christol and Mebkhout use by analogy with the theory of
classical differential modules (in which Newton polygons of
certain $t$-adic polynomials play a role), to avoid confusion
with the ``Frobenius slopes'' that also inhabit the theory
of $p$-adic differential equations, as in \cite{me-local}. 
Instead, we use the term ``break'' to evoke the connection with
ramification breaks, as encapsulated in Theorem~\ref{T:breaks}.
\end{remark}

To relate breaks to monodromy, we start in the unipotent case.
\begin{lemma} \label{L:cmunip}
For $\calE \in \calM$ unipotent, $\calE$ has highest break $0$.
\end{lemma}
\begin{proof}
By Proposition~\ref{P:characterize unip},
we can find a basis $\bv_1, \dots, \bv_n$ of $\calE$ such that
$t \frac{d}{dt} \bv_i = \sum_{j<i} c_{ij} \bv_j$
with $c_{ij} \in K$.
Since
\[
\frac{dt}{t} = \frac{d(t - t_\rho)}{t_\rho + (t - t_\rho)}
= \frac{d(t-t_\rho)}{t_\rho} \sum_{i=0}^\infty \left( 
\frac{t_\rho - t}{t_\rho} \right)^i,
\]
$\calE$ has a full basis of solutions on any disc of the form
$|t-t_\rho| < \rho$, where $|t_\rho| = \rho$. This proves the claim.
\end{proof}

\begin{lemma} \label{L:timesunip}
For $\calE, \calF \in \calM$ and $\alpha, \beta \in \RR_{\geq 0}$
with $\alpha < \beta$,
if $\calE$ has highest break $\leq \alpha$ and
$\calF$ has highest break $\beta$, then
$\calE \otimes \calF$ also has highest break $\beta$.
\end{lemma}
\begin{proof}
Let $\bv_1, \dots, \bv_m$ and $\bw_1, \dots, \bw_n$ be bases of local
horizontal sections of $\calE$ and $\calF$, respectively, around a
generic point $t_\rho$.
Then a basis of local horizontal sections of $\calE \otimes \calF$ around $t_\rho$
is given by $\bv_i \otimes \bw_j$ for $i=1,\dots, m$ and $j=1, \dots, n$.
By Lemma~\ref{L:cmunip}, we have that
$\bv_1, \dots, \bv_m$ converge for
$|t-t_\rho| < \rho^{\alpha+1}$.
Hence for $r \geq \alpha+1$,
all of the $\bw_j$ converge for $|t-t_\rho| < \rho^r$ if and only if
all of the $\bv_i \otimes \bw_j$ converge there.
This yields the desired result.
\end{proof}

We next observe how breaks are affected by a tame extension of
$k((t))$.
\begin{lemma} \label{L:tame}
Let $f: A(0,1) \to A(0,1)$ be the formally \'etale cover defined by
$f^*(t) = t^n$, for $n$ a positive integer not divisible by $p$.
For $\calE \in \calM$ and $\beta \in \RR_{\geq 0}$,
$\calE$ has highest break $\geq \beta$ (resp.\ $\leq \beta$) if and only if
$f^* \calE$ has highest break $n \beta$ (resp.\ $\leq n\beta$) if and only if
$f_* \calE$ has highest break $\beta/n$ (resp.\ $\leq \beta/n$).
\end{lemma}
\begin{proof}
This follows at once from Lemma~\ref{L:generic tame}.
\end{proof}

\begin{lemma} \label{L:wild1}
Given $a \in k^*$ and $n$ a positive integer not divisible by $p$,
let $f: X \to A[\rho_0,1)$ be a formally \'etale cover with 
reduction $k((t))[u]/(u^p - u - at^{-n})$.
Suppose $\calE \in \calM$ has no trivial submodules, but
$f^* \calE$ is constant.
Then $\calE$ has highest break $n$.
\end{lemma}
\begin{proof}
We may assume without loss of generality that 
$K$ contains an element $\pi$ with $\pi^{p-1} = -p$,
and that $\calE$ is irreducible and nonconstant. 
Then $\calE$ has rank 1, and is equal to one of the nontrivial
summands of $f_* f^* \calO$. 
If we write
\[
\Gamma(\calO,X) = \calR_{[\rho_0,1)}[u]/(1 + p \pi b t^{-n} - (1 + \pi u)^p)
\]
for $b \in \gotho_K$ reducing to $a$ in $k$,
then the trivial summand is generated by $1$ and the
nontrivial summands are generated by $(1 + \pi u)^i$
for $i=1, \dots, p-1$.
Since
\[
p \frac{d(1 + \pi u)^i}{(1 + \pi u)^i} = \frac{d(1 + 
p \pi b t^{-n})}{1 + p \pi b t^{-n}},
\]
each summand is isomorphic to the rank one $\nabla$-module generated by
$\bv$ with
\[
\nabla \bv = - \pi nb t^{-n-1} (1 + p \pi b t^{-n})^{-1} \bv \otimes dt.
\]
However, since $\exp(x)$ converges for $|x| < |\pi|$, this $\nabla$-module
is isomorphic to the $\nabla$-module generated by $\bw$ with
\[
\nabla \bw = -\pi nb t^{-n-1} \bw \otimes dt
\]
over a suitable annulus (namely $A[\rho,1)$ for $\rho = \max\{\rho_0,
p^{-1/n}\}$).

When we expand around a generic point $t_{\rho}$ with 
$|t_\rho| = \rho$, we get a local horizontal
section given by
\[
\exp ( - \pi b(t^{-n} - t_\rho^{-n}))\bw.
\]
This section converges in the disc $|t - t_\rho| < \rho^{\beta +1}$ if and only if
$|t^{-n} - t_\rho^{-n}| < 1$ throughout the disc.
However, by Lemma~\ref{L:generic tame},
the discs $|t^{-n} - t_{\rho}^{-n}|< 1$ and
$|t - t_\rho| < \rho^{n+1}$ are identical.
Thus $\calE$ has highest break $n$, as desired.
\end{proof}

\begin{lemma} \label{L:wild2}
Given $a \in k^*$ and $n$ a positive integer not divisible by $p$,
let $f: X \to A[\rho_0,1)$ be a formally \'etale cover with 
reduction $k((t))[u]/(u^p - u - at^{-n})$.
For $\beta \geq n$, 
$\calE$ has highest break $\leq \beta$ if and only if
$f^* \calE$ has highest break $\leq \beta p - n (p-1)$.
\end{lemma}
\begin{proof}
There is no loss of generality in enlarging $k$ by adjoining an
$n$-th root of $a$, and so by changing series parameters, we
may reduce to the case $a=1$. Also, we will not use the precise value
of $\rho_0$, so there is no harm in allowing it to increase
(by stipulating the existence of certain formally \'etale covers over
$A[\rho_0,1)$).

In case $n=1$, the claim follows at once from Lemma~\ref{L:generic wild}.
In the general case, let $g: A[\rho_0,1) \to A[\rho_0^n,1)$ be the $n$-th power
map. Let $f_0: X' \to A[\rho_0,1)$ be a formally \'etale cover
with reduction $k((t))[u]/(u^p - u - t^{-1})$;
then $f \circ g$ factors as
$g_0 \circ f_0$, where $g_0: X \to X'$ is a formally \'etale cover
whose reduction is totally tamely ramified of degree $n$. 

By the $n=1$
case together with Lemma~\ref{L:tame}, the following assertions are equivalent.
\begin{itemize}
\item $\calE$ has highest break $\leq \beta$.
\item $g_* \calE$ has highest break $\leq \beta/n$.
\item $f_0^* g_* \calE$ has highest break $\leq \beta p/n - (p-1)$.
\item $g_0^* f_0^* g_* \calE \cong f^* g^* g_* \calE$ has highest break
$\leq \beta p - n(p-1)$.
\end{itemize}
Moreover, there are maps
\[
\calE \to g^* g_* \calE \to \calE
\]
whose composition is multiplication by $n$,
given by the trace and adjunction maps (for a finite \'etale ring extension);
in particular, $\calE$ is isomorphic to a subobject of $g^* g_* \calE$.
Thus if $\calE$ has highest break $\leq \beta$, then
$f^* \calE$ has highest break $\leq \beta p - n(p-1)$.

On the other hand, the following assertions are also equivalent.
\begin{itemize}
\item $f^* \calE$ has highest break $\leq \beta p - n(p-1)$.
\item $(g_0)_* f^* \calE$ has highest break $\leq \beta p/n - (p-1)$.
\item $(f_0)_* (g_0)_* f^* \calE \cong g_* f_* f^* \calE$ has highest
break $\leq \beta/n$. (We did not directly address $(f_0)_*$ above, but we may
apply Lemma~\ref{L:generic wild} to deduce this.)
\item $g^* g_* f_* f^* \calE$ has highest
break $\leq \beta$.  
\end{itemize}
Moreover, there are maps
\[
\calE \to f_* f^* \calE \to \calE
\]
whose composition is multiplication by $p$,
given by the adjunction and trace maps (for a finite \'etale ring extension);
in particular, $\calE$ is isomorphic to a subobject of $f_* f^* \calE$,
and likewise of $g^* g_* f_* f^* \calE$.
Thus if $f^* \calE$ has highest break $\leq \beta p - n(p-1)$, then
$\calE$ has highest break $\leq \beta$.
\end{proof}

\begin{lemma} \label{L:wild3}
Given a polynomial $P(t^{-1})$ over $k$ of degree $n$ coprime to $p$,
let $f: X \to A[\rho_0,1)$ be a formally \'etale cover with 
reduction $k((t))[u]/(u^p - u - P(t^{-1}))$.
Suppose that $\calE \in \calM_K$ has no trivial submodules, but
$f^* \calE$ is constant.
Then $\calE$ has highest break $n$.
\end{lemma}
\begin{proof}
As in Lemma~\ref{L:wild1}, we may assume that $K$ contains
$\pi$ such that $\pi^{p-1} = -p$, and that $\calE$ is 
irreducible and nonconstant; in particular, $\calE$ is of rank 1. 
We may also assume $k$ is algebraically
closed, and that $P(t^{-1}) = \sum c_i t^{-i}$ where $c_i = 0$ whenever $i$
is divisible by $p$.

We proceed by induction on $n$. Write $P(t^{-1}) = a_n t^{-n} + Q(t^{-1})$,
where $Q$ has degree $d < n$ coprime to $p$.
Then we can write $\calE = \calF \otimes \calG$, where $\calF$
is nonconstant but
becomes constant over $k((t))[y]/(y^p - y - Q(t^{-1}))$,
and $\calG$ is nonconstant but
becomes constant over $k((t))[z]/(z^p - z - a_n t^{-n})$.
By the induction hypothesis, $\calF$ has highest break $d$;
by Lemma~\ref{L:wild1}, $\calG$ has highest break $n$.
By Lemma~\ref{L:timesunip}, $\calE$ also has highest break $n$,
as desired.
\end{proof}

At last, we can relate generic radii of convergence to ramification filtrations.
\begin{theorem}\label{T:breaks}
Suppose that $k$ is perfect.
For $\calE \in \calM^{\qu}$ whose monodromy representation has highest break
$\beta$, $\calE$ has highest break $\beta$.
\end{theorem}
\begin{proof}
Assume without loss of generality that $k$ is algebraically closed.
Let $\tau: G \to \GL(V)$ be the corresponding representation.
By Lemma~\ref{L:tame}, we may reduce to the case where
the image $H$ of $\tau$ is a $p$-group. Moreover,
thanks to Theorem~\ref{T:Matsuda}, Lemma~\ref{L:cmunip},
and Lemma~\ref{L:timesunip},
we may assume that $\calE$ is quasi-constant.

We induct primarily on the order of $H$;
the base case of order $p$ is handled by
Lemma~\ref{L:wild3}. We induct secondarily
on the lowest break $i$ of
the ramification filtration induced on $H$. 
Note that by Proposition~\ref{P:Hasse-Arf}, $i$ is an integer.

Assume without loss of generality that $\calE$ is irreducible.
Choose a maximal subgroup $N$ of $H$ containing $H^i$.
Then $N$ is normal of index $p$,
so it fixes an Artin-Schreier extension
of $k((t))$, say $k((t))[u]/(u^p - u - f)$. We can take $f$
of the form $\sum_{j=1}^i c_j t^{-j}$, with $c_j \in k$
and $c_j=0$ when $j$ is divisible by $p$. (We may omit the term
$c_0$ because $k$ is algebraically closed.)

We now pass from $k((t))$ to $E = k((t))[x]/(x^p - x - c_i t^{-i})$.
If $\calE$ becomes constant, then
Lemma~\ref{L:wild1} implies that $\calE$ has highest break
$i$, completing the induction in this case. Otherwise, let $\beta$
be the highest ramification break over $k((t))$;
then the highest ramification break over $E$ is $\beta p - (p-1)i$
(as in Example~\ref{E:as}).
If $\tau$ restricted to $\Gal(E^{\sep}/E)$
has image $H$, its lowest
break must be strictly less than $i$; otherwise, the image of $\tau$
is a proper subgroup of $H$.
Hence the induction hypothesis applies, so that
$\calE$ has highest break $\beta p - (p-1)i$ over $E$.
By Lemma~\ref{L:wild2}, $\calE$ has highest break $\leq \beta$,
with equality as long as $\beta \neq i$. But we cannot have
$\beta = i$: otherwise $H$ would be elementary abelian, hence of order
$p$ since $\calE$ is irreducible, but that case is the base case
which we handled with Lemma~\ref{L:wild3}. Thus $\beta > i$,
and $\calE$ indeed has
highest break $\beta$.
This completes both of the inductions and yields the
claim.
\end{proof}
\begin{cor} \label{C:breaks}
The equivalence of categories given by Theorem~\ref{T:monodromy rep}
commutes with the formation of break decompositions.
\end{cor}

\begin{remark}
Given Theorem~\ref{T:breaks}, 
the filtration $\calE_i$ of $\calE$ corresponding to the
break filtration has the following property:
the rank of $\calE_i$ is,
for $\rho \in (0,1)$ sufficiently close to 1, the
maximum number of linearly independent horizontal sections of $\calE$
on the open disc of radius $\rho^{i+1}$ around a generic point $t_\rho$.
\end{remark}

\begin{remark} \label{R:breaks}
In the case of $K$ discretely valued,
Theorem~\ref{T:breaks} (stated in the equivalent form of 
Corollary~\ref{C:breaks}) is due
to Matsuda \cite[Corollary~8.8]{matsuda-katz}.
Substantively similar results (comparing the ``irregularity'' of the
differential equation with the Swan conductor of the monodromy
representation) have been given by 
Tsuzuki \cite[Theorem~7.2.2]{tsuzuki-swan} and
Crew \cite[Theorem~5.4]{crew-swan}.
For more discussion, as well as a 
Tannakian reformulation, see
\cite[Complement~7.1.2]{andre} and subsequent remarks.
\end{remark}

\begin{remark}
Note that our approach to Theorem~\ref{T:breaks} is highly revisionist;
the original construction by Christol and Mebkhout of the break
decomposition of an overconvergent $\nabla$-module is by topological
means, inspired by classical techniques for studying ordinary
(complex) differential equations. 
We do not know whether our derivation of
Theorem~\ref{T:breaks} was known to Christol and Mebkhout,
or if so, whether it motivated their construction.
However, our argument by ``breaking down the module''
with successive small extensions
is loosely styled after
the derivation of the $p$-adic local monodromy theorem
given by Mebkhout,
specifically the proof of \cite[Th\'eor\`eme~5.0-20]{mebkhout}.
\end{remark}

\section{Frobenius structures on differential equations}
\label{sec:frob}

In this chapter, we recall the notion of a Frobenius structure on
a differential equation, and explain how it interacts with some of the
other structures we have introduced.

\subsection{Frobenius structures}

Recall the notion of a Frobenius structure on a $\nabla$-module.
\begin{defn}
Let $\sigma_K: K \to K$ be a continuous homomorphism acting modulo
$\gothm_K$ as the $q$-th power Frobenius, for some power $q = p^n$ of $p$.
A \emph{Frobenius lift (of order $n$)} extending $\sigma_K$ is a map $\sigma:
A[a,1) \to A[a^q,1)$ for some $a \in (0,1) \cap \Gamma^*$
which induces a map on $\calR^{\inte}$ which reduces modulo
$\gothm_{\calR}$ to the $q$-th power map.
Note that a Frobenius lift pulls back to a Frobenius lift on
any formally \'etale cover.
\end{defn}

\begin{defn} \label{D:Frobenius structure}
Given a Frobenius lift $\sigma$ of order $n$, a \emph{Frobenius structure
(of order $n$)}
on $\calE \in \calM$ is an isomorphism $F: \sigma^* \calE \to \calE$;
we will often view such an $F$ as a $\sigma$-linear map on
sections of $\calE$.
It will follow from Lemma~\ref{L:Frobenius convergent} below
(and Remark~\ref{L:convergent pullback}) that the choice of
a Frobenius structure for a single $\sigma$ determines a Frobenius
structure for every $\sigma$.
\end{defn}

\begin{lemma} \label{L:Frob bases}
Suppose $\calE \in \calM$ admits a Frobenius structure.
Then for any $\lambda \in (0,1)$, there exists $\rho_0 \in (0,1)$ such
that for any $\rho \in (\rho_0,1)$, there exists a closed aligned
subinterval $I$ of $(\rho_0,1)$ containing $\rho$
and a basis $\be_1, \dots, \be_n$ of
$\Gamma(\calE, A(I))$, such that the matrix $N$ 
over $\calR_I$ defined by
\[
t \frac{d}{dt} \be_j = \sum_i N_{ij} \be_i
\]
satisfies $|N|_r < \lambda$ for $r \in I$.
\end{lemma}
\begin{proof}
Fix a choice of a $\nabla$-module representing $\calE$, which
we will hereafter confound with $\calE$.
Choose a closed aligned subinterval $I_0 = [a,b]$ of $[0,1)$ on which $\calE$
is defined, such that $b > a^{1/q}$, and $\sigma$
maps $A[a,1) \to A[a^q,1)$.
For $l=0,1,\dots$, put $I_l = [a^{q^{-l}}, b^{q^{-l}}]$;
given a basis $\be_1, \dots, \be_n$ of
$\Gamma(\calE, A(I_0))$, define the matrix $N_l$ by
\[
t \frac{d}{dt} F^l \be_j = \sum_i (N_l)_{ij} F^l\be_i.
\]
Then
\[
N_{l+1} = \frac{dt^\sigma/dt}{t^\sigma/t} N_l^\sigma.
\]
Put $u = \frac{dt^\sigma/dt}{t^\sigma/t}$; 
then there exists $c \in (0,1)$ such that
for $\rho$ sufficiently close to 1, we have
$|N_l^\sigma|_{\rho^{1/q}} = |N_l|_\rho$
and $|u|_\rho \leq c$. It follows that
for $l$ sufficiently large, $|N_l|_\rho < \lambda$
for $\rho \in I_l$. Since the intervals $I_l, I_{l+1}, \dots$
overlap (because $b > a^{1/q}$), they cover $[a^{q^{-l}}, 1)$;
thus we have the desired result.
\end{proof}

As in \cite[Th\'eor\`eme~2.5.7]{ber2}, one has the following.
\begin{lemma} \label{L:Frobenius convergent}
If $\calE \in \calM$ admits a Frobenius structure, then
$\calE$ is overconvergent.
\end{lemma}
\begin{proof}
We need to show that for any $\eta \in (0,1)$, 
there exists $\rho \in [0,1)$ such that for any
closed aligned subinterval $I$ of $[\rho,1)$ and any
$\bv \in \Gamma(\calE, A(I))$, the sequence
\[
\left\{ \frac{1}{j!} \frac{d^j}{dt^j} \bv \right\}_{j=0}^\infty
\]
is $\eta$-null. We can rewrite the $j$-th term of this sequence as
\[
\frac{t^{-j}}{j!} \prod_{i=0}^{j-1} \left( t \frac{d}{dt} - i \right) \bv. 
\]
Suppose that $\Gamma(\calE, A(I))$ admits a basis on which
$t \frac{d}{dt}$ acts via a matrix $N$ with $|N|_r \leq |(p^m)!|$ for
some integer $m$. Then for any integer $n$,
\[
\prod_{i=0}^{j-1} \left( t \frac{d}{dt} - (np^m + i) \right)
\]
acts on this basis via a matrix $N_n$ with $|N_n|_r \leq |(p^m)!|$.
On such a basis, $\frac{t^{-j}}{j!} \prod_{i=0}^{j-1} \left( t \frac{d}{dt} - i \right)$
acts via a matrix $N_j$ with
\begin{align*}
|N_j|_r &\leq r^{-j} \frac{|(p^m)!|^{\lfloor j/p^m \rfloor}}{|j!|} \\
&\leq r^{-j} p^{-(j/p^m)((p^m-1)/(p-1)) + j/(p-1)} \\
&= r^{-j} p^{j/(p^m(p-1))}.
\end{align*}

Choose $m$ such that $p^{-1/(p^m(p-1))} > \eta$.
By Lemma~\ref{L:Frob bases}, we can find $\rho \in (0,1)$ such that
we can cover
$[\rho,1)$ with closed aligned subintervals $I$ such that each
$\Gamma(\calE, A(I))$ admits a basis on which
$t \frac{d}{dt}$ acts via a matrix $N$ with $|N|_r \leq |(p^m)!|$ for
each $r \in I$.
By increasing $\rho$ if needed, we can ensure that $\rho^{-1} 
p^{1/(p^m(p-1))} \eta < 1$. Then the matrices $N_j$ form an
$\eta$-null sequence;
it follows that $\calE$ is $\eta$-convergent, as desired.
\end{proof}

The presence of a Frobenius structure provides some restriction
on the monodromy representation, as follows.
\begin{prop} \label{P:bound tame}
Given the residual characteristic $p$ of $k$, and positive
integers $n$ and $d$,
there exists an integer $N = N(n,p,d)$ with the following property.
Suppose that $\calE \in \calM^{\qu}$ has rank $d$, and
that $\calE$ admits a Frobenius structure $F$ of order $n$. 
Then for any formally \'etale cover
$f: X \to A[a,1)$ over which $\calE$ becomes unipotent,
with reduction $R$,
the prime-to-$p$ order of the image of the inertia
subgroup of $\Gal(R/k((t)))$ in
$\Aut(H^0(\calE,X))$ is at most $N$.
\end{prop}
\begin{proof}
Note that the canonical minimal filtration of $\calE$ (in which
the first step is the maximal quasi-constant submodule, and so on) 
is preserved
by any Frobenius structure, so it suffices to consider the case
where $\calE$ is quasi-constant, then apply that case to each
successive quotient. Also, we may 
assume without loss of generality that $k$ is algebraically closed, since
the desired conclusion is insensitive to changing $K$.

Put $G = \Gal(R/k((t)))$ and $V = H^0(\calE,X)$; 
there is no harm in assuming that
$G$ injects into $\Aut(V)$. By Jordan's theorem 
(Proposition~\ref{P:Jordan}),
we can find a commutative normal subgroup $H$ of $G$ such that
$|G/H|$ is bounded as a function of $d$ alone.
By passing to the fixed field of $H$, we may reduce to the case 
where $G$ is abelian. We may also enlarge $K$ so that
all characters of $G$ are defined over $K$.

Let $\chi_1, \dots, \chi_d$ be the characters via which $G$ acts on
$V$. The fact that $F$ commutes with $G$ means that the map
$\chi \mapsto \chi^\sigma$ permutes the $\chi_i$; in particular,
$\chi_i^{\sigma^{d!}} = \chi_i$ for each $i$. It follows that
$\chi_i$ takes values in the fixed field of $K$ under $\sigma_K^{d!}$;
the group of prime-to-$p$ roots of unity in that fixed field is 
isomorphic to the group of prime-to-$p$ roots of unity in 
$\FF_{q^{d!}}$ via the Teichm\"uller map. In particular, each
$\chi_i$ has prime-to-$p$ order dividing $q^{d!} - 1$, so the 
prime-to-$p$ order of $G$ is bounded by $(q^{d!}-1)^d$, as desired.
\end{proof}
\begin{remark}
In case $K$ is discretely valued, it can be shown (using the fact
that the kernel of $\GL_d(\gotho_K) \to \GL_d(\gotho_K/\gothm_K^i)$
is torsion-free if the exponential map converges on $\gothm_K^i$)
that the entire order of the image of 
inertia is bounded by a function of $p,n,d$
and the absolute ramification index $e$ of $K$.
Indeed, that observation is implicit in Tsuzuki's proof of the unit-root
case of the $p$-adic local monodromy theorem
\cite{tsuzuki-unitroot}.
Alternatively, one can bound the $p$-part of the order of inertia by
bounding the ramification breaks of the monodromy representation
using the Christol-Mebkhout construction, then applying
Proposition~\ref{P:bound image}.
\end{remark}

\begin{remark}
Note that one cannot obtain any bound in the spirit of
Proposition~\ref{P:bound tame} in the absence of a Frobenius structure.
For instance, the rank one $\nabla$-module given by $\bv \mapsto r \bv \otimes \frac{dt}{t}$
is quasi-constant for any $r \in \QQ$ whose denominator $d$ is coprime to $r$,
but the cover required has degree $d$.
\end{remark}

\subsection{The $p$-adic local monodromy theorem}

We now give 
a form of the $p$-adic local monodromy theorem ($p$LMT) of
Andr\'e \cite{andre}, Mebkhout \cite{mebkhout}, and the present
author \cite{me-local}, which gives a context in which our results
concerning quasi-unipotent $\nabla$-modules
can be applied, at least when $K$ is discretely
valued. Our ``proof'' is nothing more than a derivation of our particular
statement from the form of the $p$LMT given in \cite{me-local}.
We follow up with a series of remarks on potential variant
forms of the $p$LMT.
\begin{theorem}[$p$-adic local monodromy theorem]
Suppose that $K$ is discretely valued. Then any 
$\calE \in \calM_K$ which admits a Frobenius structure is quasi-unipotent.
\end{theorem}
\begin{proof}
As per Definition~\ref{D:Frobenius structure}, we may assume the Frobenius
structure is with respect to any prescribed 
$\sigma$. In particular, we may assume
that $\sigma$ is a power of a Frobenius lift of order 1 (so that
we may apply the results of \cite{me-local}).
For $n$ a nonnegative integer, put $K^n = K^{\sigma_K^{-n}}$; that is,
$K^n$ is a copy of $K$ viewed as a $K$-algebra via $\sigma_K^n$.
Then by \cite[Theorem~6.12]{me-local}, there exists a
formally \'etale cover $f^n: X^n \to A_{K^n}[a,1)$ for some $n$,
and some $a \in [0,1)$, such that $\calE$ becomes unipotent on $X$.

The image of $\Gamma(\calO,X^n)$ under $\sigma^n$ induces a
formally \'etale cover $f: X \to A_K[a,1)$ whose pullback along
$A_{K^n}[a,1) \to A_K[a,1)$ coincides with $f^n$, and the images
of horizontal elements of $\Gamma((f^n)^* \calE, X^n)$ under
the $n$-th power of the Frobenius structure are horizontal
elements of $\Gamma(f^*\calE,X)$. Thus $f^* \calE$ admits
a nontrivial constant submodule; quotienting and repeating,
we deduce that $f^* \calE$ is actually unipotent over $X$.
This completes the argument.
\end{proof}
\begin{remark}
Note that the references cited above
all assume in some fashion that $k$ is perfect
(or even algebraically closed); this is so that, in our terminology,
the cover $X$ can be identified with an annulus over a finite
unramified extension of $K$. While this sort of identification
is convenient for making intermediate calculations, it does not
intervene in the definition of unipotence that we are using.
For a more robust treatment of the $p$LMT (following but simplifying
the approach of \cite{me-local}), see \cite{me-slope}.
\end{remark}

\begin{remark}
One may wonder whether the $p$-adic local monodromy theorem holds for
$K$ which are not discretely valued. Proving the $p$LMT over an
extension of $K$ proves the $p$LMT over $K$
(e.g., by the argument of \cite{me-local}*{Proposition~6.11};
see also \cite[Chapter~3]{me-part1}),
so it is enough to consider
spherically complete $K$. 
The discussion of this point naturally bifurcates following the
two proof strategies available for proving the $p$LMT;
see the two subsequent remarks.
\end{remark}

\begin{remark}
The generalization of the $p$LMT to spherically complete $K$ via the
``differential'' approach used by Andr\'e \cite{andre} and 
Mebkhout \cite{mebkhout} can probably be generalized without much
additional effort.
This approach is based on the Christol-Mebkhout
$p$-adic index theorem \cite[Section~6]{cm4}, which is known to
hold for spherically complete $K$. Although neither Andr\'e
nor Mebkhout asserts the $p$LMT for spherically complete
$K$, it seems not so difficult to check that their arguments carry
over; a key point in Andr\'e's argument is the fact that Crew
\cite[Theorem~4.11]{crew-rank1} verified the $p$LMT in rank 1,
and this argument carries over essentially unchanged.
\end{remark}

\begin{remark}
The generalization of the $p$LMT to spherically complete $K$ via the
``Frobenius'' approach used by the present author \cite{me-local}
can probably be generalized, but a fair bit of
additional effort will be required.
For one thing, it depends on Tsuzuki's unit-root $p$-adic local monodromy 
theorem \cite[Theorem~4.2.6]{tsuzuki-unitroot}, whose proof has only
partly been generalized to arbitrary $K$ by Christol
\cite{christol}; the argument founders on a thorny technical
problem that Christol did not see how to resolve, and neither do we
\cite[Remark~14]{christol}.
For another thing, the intermediate structural results of
\cite{me-local} would have to be generalized; in some cases this
seems tractable (e.g., the B\'ezout property for ``analytic rings''
can probably be generalized by imitating the corresponding arguments
of Lazard \cite{lazard}, and the ``existence of eigenvectors'' should
be straightforward), but in some cases it is less clear how to
avoid leaning on the discreteness hypothesis
(e.g., when ``raising the Newton polygon''). This issue is discussed
in somewhat more detail throughout \cite{me-slope}.
\end{remark}

\begin{remark}
It is conceivable that a version of the $p$LMT holds for overconvergent
$\nabla$-modules not necessarily admitting a Frobenius structure;
its conclusion would assert (by loose analogy with
Jordan's theorem, a/k/a Proposition~\ref{P:Jordan})
that after pulling back along a
suitable formally \'etale cover, the given $\nabla$-module
splits as a direct sum of rank 1 $\nabla$-modules, each of the form
$\nabla \bv = c \bv \otimes \frac{dt}{t}$ for some $c \in K$.
Under an additional technical hypothesis (the property ``NLE'' in
the terminology of \cite{mebkhout}), and still assuming that $K$ is
discretely valued (though perhaps the argument can be extended to
allow $K$ to be spherically complete), this is essentially the
``th\'eor\`eme de Turrittin'' of Mebkhout 
\cite[Th\'eor\`eme-5.0-20]{mebkhout}; the technical hypothesis
forces the numbers $c$ to be $p$-adically non-Liouville.
Unfortunately, it seems the only way to verify this hypothesis in
practice is by exhibiting a Frobenius structure.
Moreover, even if the Turrittin theorem held without the
non-Liouville hypothesis, it is not clear that one would be able to
avoid the use of Frobenius structures elsewhere in the theory, e.g.,
in the ``full faithfulness of overconvergent-to-convergent
restriction'' \cite{me-full}.
%it would lead
%to a stronger form of semistable reduction (i.e., the existence
%of logarithmic extensio; that is because the
%Frobenius structure is used in several other places in the argument,
%in order to provide some key finiteness arguments. 
(Thanks to Nobuo Tsuzuki for prompting this remark.)
\end{remark}

\subsection{Frobenius antecedents}

The Frobenius antecedent theorem of Christol-Dwork 
\cite[Th\'eor\`eme~5.4]{christol-dwork} gives a sufficient
condition for a $\nabla$-module on an annulus to be a Frobenius pullback
of another $\nabla$-module. The approach in \cite{christol-dwork} is 
via explicit calculations with cyclic vectors; we will instead
use Taylor series.

\begin{convention} \label{conv:standard}
Throughout this section, fix a Frobenius lift $\sigma_K$ of order $1$ on $K$,
and let $\sigma$ denote the standard extension of $\sigma_K$ to $\calR_I$
for each interval $I$; that is,
\[
\left( \sum c_i t^i \right)^\sigma = \sum_i c_i^{\sigma_K} t^{pi}.
\]
\end{convention}

\begin{theorem} \label{T:Frobenius antecedent}
For $\rho_0 \in (0,1) \cap \Gamma^*$, let $\calE$ be a
$\nabla$-module on $A[\rho_0,1)$ with the property that
\[
R(\calE, \rho) > p^{-1/(p-1)} \rho \qquad \mbox{for all
$\rho \in [\rho_0, 1)$}.
\]
Then there exists a unique $\nabla$-module $\calF$ on
$A[\rho_0^p, 1)$ such that $\calE \cong \sigma^* \calF$
and
\[
R(\calF, \rho^p) = R(\calE, \rho)^p \qquad \mbox{for all
$\rho \in [\rho_0,1)$}.
\]
\end{theorem}
\begin{proof}
Thanks to the uniqueness assertion in the claim, it suffices
to prove the claim after adjoining to $K$ a primitive $p$-th
root of unity $\zeta$.

For $i=0, \dots, p-1$, let $g_i: A[0,1) \to A[0,1)$ be the map
induced by the ring map $t \mapsto t \zeta^i$. Then
the map $h_i: g_i^* \calE \to \calE$ defined by 
\[
h_i(\bv) = \sum_{n=0}^\infty \frac{(\zeta^i-1)^n t^n}{n!} \otimes
\frac{d^n}{dt^n} \bv
\]
is well-defined because $R(\calE, \rho) > p^{-1/(p-1)} \rho$,
and is in fact an isomorphism. (Note that this is an example of the 
functoriality of rigid cohomology, as described in \cite{ber2}; see
also \cite[Proposition~2.8.1]{me-part1}.)

We may interpret the $h_i$ as giving rise to endomorphisms
of $\calE$ which are semilinear in the sense that $h_i(t\bv) = \zeta^i
t h_i(\bv)$. As in \cite[Section~5.2]{me-part1},
we set
\[
f_i(\bv) = t^{-i} \sum_{e=0}^{p-1} \zeta^{-ei} h_e(\bv),
\]
and note that $f_i(\bv)$ is fixed by the $h_i$; we also note that
if we apply each $f_i$ to each of a set of generators of $\calE$,
the images generate $\calE$ over $\calO^\sigma$. In other words,
the images generate a $\nabla$-module $\calF$ over 
$A[\rho_0^p,1)$ with $\sigma^* \calF \cong \calE$.
Moreover, if $t_\rho \in \CC$ is a generic point of radius $\rho$,
then we can apply the $f_i$ to horizontal sections of
$\calE$ on the disc $|t - t_\rho| < \rho^{\beta + 1}$
to obtain horizontal sections of $\calF$ on the disc
$|t - t_\rho^p| < \rho^{p(\beta + 1)}$.
We thus have $R(\calF, \rho^p) \geq R(\calE,\rho)^p$; we have
the reverse inequality also by Lemma~\ref{L:generic Frob}. Thus
\[
R(\calF, \rho^p) = R(\calE,\rho)^p.
\]

To conclude, we verify uniqueness. To do this, it suffices to show that
if there exists $\calF'$ such that
$\calE \cong \sigma^* \calF'$ and $R(\calF', \rho^p) = R(\calE,\rho)^p$
for all $\rho \in [\rho_0,1) \cap \Gamma^*$, then $\calF'$ is fixed
under each of the $h_i$; once this is known, it follows that
the fixed locus of the $h_i$ is precisely $\calF'$,
and hence that $\calF \cong \calF'$.

Let $I = [a,b]$ be a closed aligned subinterval of $[\rho_0,1)$ such that
$R(\calF',\rho^p) \geq p^{-p/(p-1)} b^p$ for any $\rho \in I$,
and put $I^p = [a^p, b^p]$.
Then the Taylor series map gives an isomorphism $\pi_2^* \calF'
\to \pi_1^* \calF'$ on the subspace of $A(I^p) \times A(I^p)$,
with coordinates $t_1$ and $t_2$, where $|t_1 - t_2| \leq 
p^{-p/(p-1)} b^p$.
This isomorphism pulls back to the Taylor series isomorphism
$\pi_2^* \calE \to \pi_1^* \calE$ on the subspace of $A(I) \times
A(I)$ where $|t_1^p - t_2^p| \leq p^{-p/(p-1)} b^p$,
and in particular to the subspace where $|t_1 - t_2| \leq p^{-1/(p-1)} b$.
That means we can reconstruct the $h_i$ starting with the trivial action
on $\calF'$.
As noted above, this yields that $\calF \cong \calF'$, yielding
the desired uniqueness.
\end{proof}

\begin{remark}
Note that \cite[Th\'eor\`eme~5.4]{christol-dwork} is slightly
weaker, as it requires by hypothesis the stronger bound
$R(\calE,\rho) > p^{-1/p} \rho$. By contrast, the
condition $R(\calF, \rho^p) = R(\calE,\rho)^p$ is quite necessary,
as evidenced by the following example
\cite[5.1]{christol-dwork}:
if $\calF$ is the rank 1 $\nabla$-module given by
$\nabla \bv = (pt)^{-1} \otimes dt$, then
$\calF \not\cong \calO$ but $\sigma^* \calF \cong \sigma^* \calO$.
\end{remark}

\subsection{Frobenius and ramification}

As we have seen, the Christol-Mebkhout characterization of breaks of
a quasi-unipotent $\nabla$-module involves inspecting the rate of
convergence of horizontal sections near the boundary of a disc.
In the presence of a Frobenius structure, its
global convergence properties can be used instead.

\begin{prop} \label{P:frobbreak}
Retain Convention~\ref{conv:standard}.
Suppose that for some $\epsilon \in (0,1) \cap \Gamma^*$,
$\calE$ is a $\nabla$-module over $A[\epsilon,1)$
equipped with a Frobenius structure given by an
isomorphism $F: \sigma^* \calE \to \calE$ over
$A[\epsilon^{1/p},1)$.
Suppose further that for any closed aligned subinterval $I$ of
$[\epsilon,1)$, there exists an $\gotho$-lattice in $\Gamma(\calE, A(I))$
stable under $t \frac{\del}{\del t}$.
Then the highest break $\beta$ of $\calE$ (in the sense of 
Theorem~\ref{T:breaks}) satisfies the inequality 
\[
\beta \leq \frac{1}{(p-1) \log_p (\epsilon^{-1})}.
\]
\end{prop}
\begin{proof}
For $\rho \in (\epsilon,1) \cap \Gamma^*$,
let $t_\rho$ be a generic point of radius $\rho$.
By hypothesis, over the disc $|t - t_\rho| < \rho$, $\calE$ admits
an $\gotho$-lattice stable under $t \frac{\del}{\del t}$. Since
$(t-t_\rho)/t$ has norm less than 1 throughout this disc,
an $\gotho$-lattice stable under $t \frac{\del}{\del t}$
is also stable under $(t-t_\rho) \frac{\del}{\del(t-t_\rho)}$.
By a direct calculation,
we deduce that
\[
R(\calE, \rho) \geq p^{-1/(p-1)} \rho.
\]
By applying Frobenius (and Lemma~\ref{L:generic Frob}), 
we have $R(\calE, \rho^{1/p^m}) \geq p^{-1/(p^m(p-1))} \rho^{1/p^m}$.
If $\beta$ is the highest break of $\calE$, for large $m$ (hence for any $m$)
one then has the inequality
\[
\rho^{\beta/p^m} \geq p^{-1/(p^m(p-1))};
\]
by taking limits, we obtain the same inequality with $\rho = \epsilon$.
This yields the desired result.
\end{proof}

%\begin{remark}
%We note in passing that Theorem~\ref{T:Frobenius antecedent} applies
%to $\sigma^* \calE$ for $\rho \in [\epsilon^{1/p}, 1)$, since
%\[
%R(\sigma^* \calE, \rho) \geq R(\calE, \rho^p)^{1/p}
%\geq p^{-1/(p(p-1))} \rho > p^{-1/(p-1)} \rho.
%\]
%\end{remark}

\begin{remark}
Note that Proposition~\ref{P:frobbreak} 
shows that given the $p$LMT, one can use explicit convergence
information for Frobenius
to control the extension needed to find the unipotent basis.
It would be interesting to turn this argument on its head, and use
Frobenius convergence information to give a more direct proof of the
$p$LMT. However, we have no idea how to do this, even in the unit-root
case originally treated by Tsuzuki \cite{tsuzuki-unitroot};
such an approach would suggest a method for extending Tsuzuki's
arguments to the case of $K$ spherically complete.
The Christol-Mebkhout $p$-adic index theorem from \cite{cm4}
does something analogous using connection convergence information,
but it does not say anything about tame ramification.
\end{remark}

%\begin{remark}
%Later in this series, we will be considering $\nabla$-modules
%with Frobenius on the product of an annulus and a disc. We will be
%blowing up in the center and examining the monodromy of the resulting
%$\nabla$-modules. In so doing, we will use
%Proposition~\ref{P:frobbreak} to give a \emph{uniform} bound on the
%highest breaks.
%\end{remark}

\subsection*{Acknowledgments}

Thanks to the Institute for Advanced Study and the
Universit\`a degli Studi di Padova (especially Francesco
Baldassarri and Bruno Chiarellotto) for their hospitality.
The author was partially supported by NSF grants DMS-0111298
and DMS-0400727.


\begin{thebibliography}{99}

\bibitem{andre}
Y. Andr\'e, Filtrations de type Hasse-Arf et monodromie $p$-adique,
\textit{Invent.\ Math.} \textbf{148}(2) (2002) 285--317.

\bibitem{balda-divizio}
F. Baldassarri and L. di Vizio, 
Continuity of the local $p$-adic radius of convergence for a 
connection on a Berkovich analytic space, in preparation.

\bibitem{berger-cst}
L. Berger, Repr\'esentations $p$-adiques et \'equations diff\'erentielles,
\textit{Invent.\ Math.} \textbf{148}(2) (2002) 219--284.

\bibitem{berger-weak}
L. Berger, \'Equations diff\'erentielles $p$-adiques et $(\phi, N)$-modules
filtr\'es, preprint, \texttt{arXiv: math.NT/0406601}, version of 29 Jun 2004.

\bibitem{berkovich-ihes}
V.G. Berkovich, \'Etale cohomology for non-Archimedean analytic spaces,
\textit{Publ.\ Math.\ IH\'ES} \textbf{78} (1993) 5--161.

\bibitem{berkovich-icm}
V.G. Berkovich, $p$-adic analytic spaces, in
Proceedings of the International Congress of Mathematicians,
Vol. II (Berlin, 1998), \textit{Doc.\ Math.} Extra Vol. II (1998) 141--151
(electronic).

\bibitem{ber2}
P. Berthelot, Cohomologie rigide et cohomologie rigide \`a support
    propre. Premi\`ere partie,
Pr\'epublication IRMAR 96-03, available at
    \texttt{http://www.maths.univ-rennes1.fr/\textasciitilde berthelo},
1996.

\bibitem{chiar-lestum}
B. Chiarellotto and B. le Stum,
Pentes en cohomologie rigide et $F$-isocristaux unipotents,
\textit{Manuscripta Math.} \textbf{100} (1999) 455--468.

\bibitem{christol}
G. Christol, About a Tsuzuki theorem,
\textit{$p$-adic Functional Analysis (Ioannina, 2000)},
Lecture Notes in Pure and Appl.\ Math.\ 222, (Dekker, 2001),
pp.\ 63--74.

\bibitem{christol-dwork}
G. Christol and B. Dwork,
Modules diff\'erentiels sur des couronnes,
\textit{Ann.\ Inst.\ Fourier (Grenoble)},
\textbf{44}(3) (1994) 663--701.

\bibitem{cm3}
G. Christol and Z. Mebkhout,
Sur le th\'eor\`eme de l'indice des \'equations
            diff\'erentielles $p$-adiques. III,
\textit{Ann.\ of Math.} \textbf{151}(2) (2000) 385--457.

\bibitem{cm4}
G. Christol and Z. Mebkhout,
Sur le th\'eor\`eme de l'indice des \'equations
            diff\'erentielles $p$-adiques. IV,
\textit{Invent.\ Math} \textbf{143}(3) (2001) 629--672.

\bibitem{crew-rank1}
R. Crew, $F$-isocrystals and $p$-adic representations,
\textit{Algebraic Geometry, Bowdoin, 1985 (Brunswick, Maine, 1985)},
Proc.\ Sympos.\ Pure Math.\ 46 (American Math.\ Society, 1987), pp.\
111--138.

\bibitem{crewfin}
R. Crew, Finiteness theorems for the cohomology of an overconvergent
            isocrystal on a curve,
\textit{Ann.\ Sci.\ \'Ecole Norm.\ Sup.\ (4)}
\textbf{31}(6) (1998) 717--763.

\bibitem{crew-swan}
R. Crew, Canonical extensions, irregularities, and the Swan conductor,
\textit{Math.\ Ann.} \textbf{316}(1) (2000) 19--37.

\bibitem{fresnel}
J. Fresnel and M. van der Put,
\textit{Rigid Analytic Geometry and Its Applications},
Progress in Mathematics 218 (Birkh\"auser, 2004).

\bibitem{gruson}
L. Gruson, Fibr\'es vectoriels sur un polydisque ultram\'etrique,
\textit{Ann.\ Sci.\ \'Ecole Norm.\ Sup.\ (4)} \textbf{1} (1968), 45--89.

\bibitem{isaacs}
I.M. Isaacs,
\textit{Character Theory of Finite Groups}
(Dover, 1994).

\bibitem{jacobson}
N. Jacobson, \textit{Basic Algebra II} (W.H. Freeman, 1989).

\bibitem{katz-gauss}
N.M. Katz, \textit{Gauss Sums, Kloosterman Sums, and Monodromy Groups},
Annals of Mathematics Studies 116 (Princeton University Press, 1988).

\bibitem{me-local}
K.S. Kedlaya, A $p$-adic local monodromy theorem,
\textit{Ann.\ of Math.} \textbf{160} (2004) 93--184.

\bibitem{me-full}
K.S. Kedlaya, Full faithfulness for overconvergent $F$-isocrystals,
\textit{Geometric Aspects of Dwork Theory} (de Gruyter, 2004), pp.\ 
819--835.

\bibitem{me-part1}
K.S. Kedlaya, Semistable reduction for overconvergent $F$-isocrystals, I:
    Unipotence and logarithmic extensions,
preprint, \texttt{arXiv: math.NT/0405069}, version of 9 Jan 2005.

\bibitem{me-slope}
K.S. Kedlaya, Slope filtrations revisited,
preprint, \texttt{arXiv: math.NT/0504024}, version of 19 Apr 2005.

\bibitem{lazard}
M. Lazard, Les z\'eros des fonctions analytiques d'une variable sur un
            corps valu\'e complet,
\textit{Publ.\ Math.\ IH\'ES} \textbf{14} (1962) 47--75.

\bibitem{matsuda-katz}
S. Matsuda, Katz correspondence for quasi-unipotent overconvergent
            isocrystals,
\textit{Comp.\ Math.} \textbf{134}(1) (2002) 1--34.

\bibitem{mebkhout}
Z. Mebkhout, Analogue $p$-adique du th\'eor\`eme de Turrittin et le
            th\'eor\`eme de la monodromie $p$-adique,
\textit{Invent.\ Math.} \textbf{148}(2) (2002) 319--351.

\bibitem{nagata}
M. Nagata, \textit{Local Rings},
Interscience Tracts in Pure and Applied Mathematics, No. 13
(John Wiley \& Sons, 1962).

\bibitem{serre}
J.-P. Serre, \textit{Local Fields} (translated from the French by
M.J. Greenberg), Graduate Texts in Mathematics 67
(Springer-Verlag, 1979).

\bibitem{tsuzuki-swan}
N. Tsuzuki, The local index and the Swan conductor,
\textit{Comp.\ Math.} \textbf{111}(3) (1998) 245--288.

\bibitem{tsuzuki-unitroot}
N. Tsuzuki, Finite local monodromy of overconvergent unit-root
            $F$-isocrystals on a curve,
\textit{Amer.\ J.\ Math.} \textbf{120}(6) (1998) 1165--1190.

\bibitem{vdp}
M. van der Put, Cohomology on affinoid spaces,
\textit{Comp.\ Math.} \textbf{45}(2) (1982) 165--198.

\end{thebibliography}
\end{document}